\documentclass[a4paper,11pt]{article}
\usepackage{setspace}
\usepackage{cancel}
\usepackage{amsmath}
\usepackage{amssymb}
\usepackage{soul}

\usepackage{arydshln}
\usepackage{adjustbox, lipsum}
\usepackage{enumerate}
\usepackage[ruled,vlined,linesnumbered]{algorithm2e}
\usepackage{rotating}
\usepackage[utf8]{inputenc}
\usepackage{amsfonts}
\usepackage{pstricks}
\usepackage{pdflscape}
\usepackage{graphicx}
\usepackage{soul}
\usepackage{afterpage}
\usepackage{url}

\usepackage{rotating}
\usepackage[sort]{natbib}

\usepackage{booktabs}

\usepackage{authblk}
\usepackage{multirow}
\usepackage{mathtools}
\usepackage{scalerel}

\newtheorem{prop}{\bf Proposition}[section]

\SetKwRepeat{Do}{do}{while}
\providecommand{\norm}[1]{\lVert#1\rVert}
\DeclarePairedDelimiter\abs{\lvert}{\rvert}%

\newcommand{\bb}{b}
\newcommand{\name}{\mbox{(RL-FS)}}
\newcommand{\nameM}{\mbox{(RL-FS-M)}}
\usepackage{xcolor}
\definecolor{gr}{RGB}{0,153,0}

\author[(a)]{Marta Baldomero-Naranjo\footnote{ Corresponding author}}
\author[(a)]{Luisa I. Martínez-Merino}
\author[(a)]{Antonio M. Rodríguez-Chía}

\affil[(a)]{\small{Departamento de Estadística e Investigación Operativa, Facultad de Ciencias, Universidad de Cádiz, Campus Universitario, Puerto Real, (Cádiz), Spain, marta.baldomero@uca.es, luisa.martinez@uca.es, antonio.rodriguezchia@uca.es}}
\date{}

\title{\LARGE A robust SVM-based approach with feature selection and outliers detection for classification problems}

\makeatletter
\let\oldabs\abs
\def\abs{\@ifstar{\oldabs}{\oldabs*}}
\makeatother
\begin{document}

%
\maketitle
\begin{abstract}
 This paper proposes a robust classification model, based on support vector machine (SVM), which simultaneously deals with outliers detection and feature selection.
The classifier is built considering the ramp loss margin error and it includes a budget constraint to limit the number of selected features. The search of this classifier is modeled using a mixed-integer formulation with big M parameters. Two different approaches (exact and heuristic) are proposed to solve the model. 
 The heuristic approach is validated by comparing the quality of the solutions provided by this approach with the exact approach.
In addition, the classifiers obtained with the heuristic method
are tested and compared with existing SVM-based models to demonstrate their efficiency.


\end{abstract}

\section{Introduction}
\label{sec:introduction}
Support vector machine (SVM) is a frequently used mathematical programming approach in supervised classification. It was introduced by \cite{13SVN} and \cite{Vapnikunicidad}. It has been applied to different fields such as finance \citep{ap7,apcreditscoring,ap5,AP-energyandfinance}; machine vision \citep{ap2,AP-machinevision}; biology \citep{ap9,Ap-agricultura,ApMedicine}; or sustainability \citep{Ap-medioambiente,AP-airport}. More applications can be found in \cite{CERVANTES2020}. 

Given a set  $N:=\{1,\ldots,n\}$ of individuals partitioned into two classes $\cal{Y}=\mbox{\{-1,1\}}$, each individual $i \in N$ is associated with a pair $(x_i,y_i) \in \mathbb{R}^d \times \{-1,1\}$, where $d$ is the number of features analyzed over each individual of $N$, $x_i$ contains the feature values, and $y_i$ provides the class membership, 1 or -1. The aim of SVM is to find a separating hyperplane $w\cdot x + \bb=0 $ to classify new individuals by finding the balance between maximizing the distance (margin) between the two parallel hyperplanes supporting individuals of each class and minimizing the error of misclassified data. See \cite{01bradley1998,BlancoPuertoChia} for the use of SVM model considering different norms to measure the margin.

Although classical SVM models have high predictive power in comparison with other state-of-the-art classifying methods, some drawbacks also arise from their use. As stated in \cite{Hastie}, classical SVMs are not robust against outliers and they do not perform efficiently when the analyzed dataset contains a large number of irrelevant features. In this paper we introduce a model which deals with both aspects: outliers detection and feature selection.

In this context, an outlier is an individual whose feature values are not similar from the corresponding ones of most individuals that belong to  the same class. The presence of outliers in real-life data is common and is caused by several factors which include but are not limited to: transcription errors, measurement errors, mislabelling, data processing errors, intentional instances (e.g. spam). There are many papers in the literature that use an SVM based model to address these issues, see for instance, \cite{out1,out2,out3,31Brooks2011,BlaJapPue20LN}. 

 Besides, \cite{Karami2018} also deals with outliers detection. Specifically, the mentioned paper considers a novel anomaly-based intrusion detection system in presence of long-range  independence data called benign outliers and also provides visualized information and insights to 
end users. This type of methodology is particularly suitable in the abnormal network traffic analysis which is an important research topic to protect computer networks from intruders. Some recent papers are focused on this topic, see \cite{Bao2016,KABIR2018303,KARAMI20151,KARAMI20151253,KARAMI20151262,Karami2014UtilizationOM}.

We will focus on the approach introduced by \cite{31Brooks2011}.
This specific model is known as SVM with ramp loss and its formulation for the $\ell_1$-norm case is given by
\begin{alignat}{4}
\mbox{(RL-$\ell_1$)} &\quad  &  & \mbox{min}&  \quad & 	\displaystyle \norm{w}_1 + C  \left( \sum_{i=1}^n\xi_i + 2 \sum_{i=1}^n z_i \right), & & \nonumber \\
& & &\mbox{s.t.} &  \quad &  \mbox{if } z_i=0,\;\; y_i \left(\sum_{k=1}^{d}w_k  x_{ik}+\bb \right) \geq 1- \xi_i,   &  & i \in N, \label{rl-cons} \\
& & & & & 0 \leq   \xi_i \leq 2, & \quad  & i \in N, \label{rl-l2-ec2}  \\
& & & & & z_i\in \{0,1\}, & &  i \in N, \label{rl-l2-ec3} 
\end{alignat}
where $w$- and $\bb$-variables are the hyperplane coefficients. Furthermore, $C$ is a non-negative constant, and $\norm{w}_1$ represents the $\ell_1$-norm of vector $w$.
In addition, $\xi_i$ determines penalization if the individual $i$ is  within the strip defined by the two parallel hyperplanes ($wx+b=1$ and $wx+b=-1$) and $z_i$ determines whether $i$ is a misclassified object outside this strip or not, taking the value 1 or 0, respectively. Observe that in this model, penalization in the objective function of each misclassified individual is bounded by two. This model has been analyzed in the literature, see for instance \cite{31Brooks2011, carrizosaheuristic,bonami,BlaJapJus19,BalMarRod19}.
In this paper, we propose an SVM model with ramp loss derived from the aforementioned model, which includes feature selection and it is formulated as a mixed integer linear program.  

As previously mentioned, a weakness of the classical SVM model is that its mathematical formulation does not limit the number of $w$-variables that are not zero in the optimal solution. It may result in classifiers that require information from many features in order to determine an  individual's class. In real-life, analyzing many features could imply high costs. For example, if we want to classify whether or not a patient has a disease, each feature usually represents a medical analysis and each of these has a cost in money and time. Therefore, our aim is to use the least number of medical tests getting the maximum accuracy when classifying the patient. Usually, it feels that the more tests are performed, the better the diagnosis. However, the performance of a large number of medical tests can lead to high costs for the health system or the patient, saturation of the health system, inconvenience for the patient (invasive tests, numerous visits to the clinic, etc.), and longer waiting time to obtain the diagnosis. For example, in the United States, there is a 62\% of personal bankruptcies incurred in medical bills, and 90\% of patients are in debt, as reported \cite{Lee2020,Facturas10}. Many of the aforementioned drawbacks could be avoided if only the tests that are truly meaningful were carried out. Moreover, classifiers are easier to interpret when the number of features is reduced rather than the ones that contain many features. 
  Furthermore, the knowledge of the features that have a relevant influence in a medical diagnosis is highly applicable in real-life because it can help in opening new lines of medical research as new treatments or earlier diagnosis.

In general, feature selection techniques have been divided into three different groups: filter, wrapper and embedded methods. Filter methods are based on a preliminary study of each feature's relevance and only features with  significant importance are considered for the classification method, see \cite{GuyGunNikZad06}. In contrast, wrapper methods interact with the classification method to select the set of relevant features, see  \cite{Kohavi,WrapperESWA20}. Finally, embedded methods study the feature selection and the classification simultaneously in the same model. Many recent papers in the literature propose different techniques to deal with feature selection in support vector machine. These models can be considered as embedded methods, see for example   \citet{21Maldonado2014,24Aytug2015,23Antonio, Ghaddar,Jimenez, Kunapuli,22Luisa,MalLopJimHer19,GauGorHir19, Nguyen,Cura2020, Lee2020}.

In this paper, we introduce an SVM with ramp loss limiting the number of selected features including a budget constraint. In this model we use the $\ell_1$-norm, the norm commonly used  in the literature for feature selection  in SVM based models. Adopting this norm results in classifiers where the number of selected features is lower than using $\ell_2$-norm due to its sparse property. Besides, the fact that the proposed model determines the outliers and selects the features simultaneously have some advantages in practical environments with respect to other approaches that do these processes independently.

One advantage of the proposed model is that it prevents the loss of valuable information due to the incorrect removal of elements of the  sample. In many  cases, several of the features of an individual in the sample may take anomalous values. However, it does not make sense
to exclude this individual if those features 
do not affect the classification. This is particularly striking when 
dealing with medical data, where usually a large sample is not available.

Another advantage is that the developed model allows us to identify as outliers data that  were not initially classified as such.
Indeed, an individual might not be identified as outlier for a general outlier detection procedure because most of its features take similar values to the ones of their class and only some of them take very different values. However, if after the feature selection, these features with very different values are the most relevant in the classification process, this individual should be classified as outlier. Precisely, in the proposed model, feature selection and outliers  detection  processes are carried out simultaneously.

In addition to the introduction of a new model, we develop  some strategies to improve the solution time. Furthermore, we provide a heuristic method that allows us to obtain good quality solutions very quickly. Concretely, we adapt a heuristic approach that has been previously used in different mixed integer linear models \citep[see][]{Guastaroba2012, AngManGra10,AKS}. This heuristic is known as Kernel Search (KS) and its main idea is to iteratively provide a better feasible solution to the problem by solving a sequence of restricted MILPs obtained from the original model.

The remainder of this paper is structured as follows. Section \ref{general_model} introduces the model, presents a valid inequality and develops a procedure to tighten the big M parameters of the formulation. Moreover, in this section, valid values for the big M parameters are computed. 
 Section \ref{sec:Heuristic} presents a heuristic based on the Kernel Search for solving the proposed model. Section \ref{sec:computationalExp}
    contains computational experiments carried out on real-life datasets. In this section,  we validate the heuristic algorithm analyzing the best solutions reported by this algorithm and the exact approach. We also test the efficiency of the proposed classifier on real-life datasets comparing its predictive power with other SVM-based models existing in the literature.
    Our conclusions and some future research topics are included in Section \ref{conclusions}.

\section{The model}
\label{general_model}


In this section, we introduce a model based on SVM that tries to eliminate the adverse effects of outliers using a small number of relevant features. Regarding feature selection, this model 
 includes a budget constraint in the formulation to limit the maximum number of features required for classification.
This budget constraint has also been used in \citet{21Maldonado2014} and \citet{22Luisa}. The model could be slightly modified associating a cost with each feature. Thus, the resulting budget constraint would restrict the cost of classifying a new individual. 
With respect to the adverse effects of outliers, this model 
uses the ramp loss margin error introduced by \cite{31Brooks2011} to avoid outliers influence.
Finally, this model considers the $\ell_1$-norm to measure the distances and it is formulated as a mixed integer program with conditional constraints:
\begin{alignat}{3}
  & \min&  \quad & 	\displaystyle \sum_{k=1}^{d}\eta_k + C  \left( \sum_{i=1}^n\xi_i + 2 \sum_{i=1}^n z_i \right), & & \nonumber \\
 &\mbox{s.t.} &\quad &\eqref{rl-cons}-\eqref{rl-l2-ec3},\nonumber\\
& & & - l_k v_k\leq w_k \leq u_k v_k, & \quad& k\in D,   \label{MILP1:ec1} \\
& & & -\eta_k\leq w_k \leq \eta_k, & \quad& k \in D, \label{ec:eta}  \\
& & & \sum_{k=1}^{d}v_k \leq B. & \quad   & \label{MILP1:ec3}  \\
& & & v_k \in \{0,1\}, & \quad   & k\in D,   \label{MILP1:ec2}\\
& & & \eta_k \geq 0, & \quad   & k\in D,\label{ec:eta2}
\end{alignat}  
where $\eta$-variables represent the absolute value of the hyperplane coefficients $w$ and $v$-variables are binary variables that indicate whether the associated feature is selected or not. Observe that constraint \eqref{MILP1:ec3} limits the number of selected features.

This model can be reformulated by expressing each $w_k$-variable for $k\in D$ as the difference between two 
non-negative variables $w_k^+$ and $w_k^-$.  Since $ w_k^+ + w_k^-,
$ for  $k \in D$ is part of the objective function, then $\left|w_k\right|= 
w_k^+ + w_k^-$ and, at most, one of the two variables for any $k\in D$  is non-zero in an optimal  solution, see for instance \cite{22Luisa} or \cite{BalMarRod19}. Hence, the model above can be equivalently reformulated as,
\begin{alignat}{4}
	\mbox{\name} & \quad &  & \min& \quad & 	\displaystyle \sum_{k=1}^{d}\left(w_k^+ + w_k^-\right) + C  \left( \sum_{i=1}^n\xi_i + 2 \sum_{i=1}^n z_i \right), & & \nonumber \\
	& & &\mbox{s.t.} &  \quad &\eqref{rl-l2-ec2}-\eqref{rl-l2-ec3}, \eqref{MILP1:ec3}-\eqref{MILP1:ec2},  & \quad & \nonumber \\ 
	& & & & &  \mbox{if } z_i=0,\;\; y_i \left(\sum_{k=1}^{d}(w_k^+-w_k^-) x_{ik}+\bb \right) \geq 1- \xi_i,     & \quad &i \in N, \label{ec:cond}  \\
	&&& & &  w_k^+ \leq u_k v_k, & \quad &k\in D,    \label{UB:w} \\
	&&& & &  w_k^- \leq l_k v_k, & \quad &k\in D,    \label{LB:w}\\
	&&& & &  w_k^+\geq 0, w_k^- \geq 0, & \quad &k \in D.   \label{ec:w} 
\end{alignat} 

Moreover, \mbox{\name} can be linearized  substituting \eqref{ec:cond} by: 
\begin{equation}
y_i\left( \sum_{k=1}^{d}(w_k^+ - w_k^-) x_{ik}+\bb\right) \geq 1- \xi_i -M_iz_i, \quad   i\in N. \label{ec:linealizado}
\end{equation}
We will refer to the resulting model as \mbox{\nameM}. 
The choice of an appropriate value for $M_i,$ for $i \in N$ is essential to provide efficient solution approaches for the model; i.e., $M_i$ should be big enough so that both models are equivalent, but it should be also as small as possible to provide good linear relaxation and to reduce the computational time for solving this model. 
In the following proposition, we give a result that establishes a relationship between $\xi$-variables and $z$-variables.

\begin{prop}\label{remark1}
 An optimal solution of \mbox{\name}, $({w^+}^*,{w^-}^*, \bb^*, \xi^*, z^*,v^*),$ satisfies the following condition:
		\begin{equation}
		\xi_i^* z_i^*=0, \label{cond:1} \quad i \in N. 
		\end{equation}
\end{prop}
The proof can be obtained by contradiction in an easy way, for this reason, it has been omitted. 
As a consequence of the previous result, the linearized version of condition \eqref{cond:1} given by 
	\begin{equation}\label{rest:cond1}
	\xi_i \leq 2(1 - z_i), \quad i \in N,
	\end{equation}
will be used to strengthen the formulations \mbox{\nameM}. Henceforth, we will refer to \mbox{\nameM}$+\eqref{rest:cond1}$ as \nameM, unless stated otherwise. 

In order to solve the model, it is necessary to provide valid values of big M parameters. Model \mbox{\nameM} contains two sets of big M parameters: the ones associated with the families of constraints \eqref{UB:w} and \eqref{LB:w}, and also the ones associated with the set of constraints \eqref{ec:linealizado}. The following  subsections will focus on initializing and improving these big M parameters in order to efficiently solve the model.

\subsection{Initial bounds for the big M parameters}

In \cite{31Brooks2011}, a mixed integer linear model for the classical SVM with ramp loss was introduced. In addition, in \citet{BalMarRod19} a result establishing theoretical bounds for the big M parameters appearing in the ramp loss model was proposed. In fact, this result could be adapted for $\nameM$ obtaining the following proposition,

\begin{prop}\label{prop:BigM}
Taking the values $M_i,$ for $i \in N,$ such that,
	$$M_i \geq \left(\underset{j \in N}{\max } \; \left\{ \norm{x_i-x_j}_q : y_i=y_j \right\} \right) \norm{w}_p ,$$ 
where $\norm{\cdot}_p$ represents the $\ell_p$-norm and $\|\cdot\|_q$ its dual, the problems $\name$ and $\nameM$ are equivalent. 
\end{prop} 

This proposition can be proven analogously to the one presented in \cite{BalMarRod19}. This result bounds the big M parameters by a value composed by two terms. The first term, $\underset{j \in N}{\max } \; \left\{ \norm{x_i-x_j}_q : y_i=y_j \right\}$, is easily computed for any norm, but computing the second term, $\norm{w}_p$, is not that easy. 
However, note that if we establish $p=1$ and we have an upper bound  (UB) on the optimal objective value of \nameM, then $\mbox{UB}\geq \norm{w}_1$. Consequently, initial big M parameters for the proposed model can be calculated.

The first step is to compute an upper bound (UB) on the optimal objective value. To do so, we solve the  classical SVM model using $\ell_1$-norm that will be denoted as (SVM-$\ell_1$) and formulated as follows:
\begin{alignat}{4}
\mbox{(SVM-$\ell_1$)} & \quad &  & \mbox{min}&  \quad & 	\displaystyle \sum_{k=1}^{d}\eta_k + C  \sum_{i=1}^n\xi_i , & & \nonumber \\
& & &\mbox{s.t.} &  \quad & \eqref{ec:eta},\eqref{ec:eta2} \nonumber\\
& & & & & y_i \left(\sum_{k=1}^{d}w_k x_{ik}+\bb \right) \geq 1- \xi_i,   &  & i \in N,  \\
& & & & &  \xi_i \geq 0, & \quad  & i \in N.
\end{alignat}
From its optimal solution, ($w^{\tiny{{\mbox{SVM}}}},\bb^{\tiny{{\mbox{SVM}}}}, \xi^{\tiny{{\mbox{SVM}}}}$), we build  $({\tilde{w}^+},{\tilde{w}^-}, \tilde{\bb}, \tilde{\xi}, \tilde{z},\tilde{v}),$ as follows: 
\begin{eqnarray*}
&\tilde{w}^+_k = \begin{cases}
w^{\tiny{{\mbox{SVM}}}}_k, &\text{if } w^{\tiny{{\mbox{SVM}}}}_k\geq 0,\\
0, &\text{otherwise,}
\end{cases}\text{for }k \in D,\;\;\;
& \tilde{w}^-_k = \begin{cases}
-w^{\tiny{{\mbox{SVM}}}}_k,&\text{if } w^{\tiny{{\mbox{SVM}}}}_k\leq 0,\\
0, &\text{otherwise,}
\end{cases}\text{for }k \in D,\;\;\;
\\
&\tilde{\xi}_i= \begin{cases}
\xi^{\tiny{{\mbox{SVM}}}}_i, &\text{if }\xi^{\tiny{{\mbox{SVM}}}}_i \leq 2,\\
0, &\text{otherwise,}
\end{cases}\text{for }i \in N,\;\;\;
& \tilde{z}_i = \begin{cases}
0, &\text{if } \xi^{\tiny{{\mbox{SVM}}}}_i \leq 2,\\
1, &\text{otherwise,}
\end{cases} \text{for }i \in N,
\\
&\tilde{v}_k= \begin{cases}
0, &\text{if }w^{\tiny{{\mbox{SVM}}}}_k= 0,\\
1, &\text{otherwise,}
\end{cases}\text{for }k \in D,\;\;\;
& \tilde{\bb}=\bb^{\tiny{{\mbox{SVM}}}}.
\end{eqnarray*}

If constraint \eqref{MILP1:ec3} is fulfilled, a feasible solution of problem $\nameM$ is obtained. If not,  we sort the components of $w^{\tiny{{\mbox{SVM}}}}$-vector in non-decreasing order and we fix to zero the first $n-B$ $w$-variables. Fixing these variables to zero, we solve the model again (SVM-$\ell_1$). From its optimal solution, we build a feasible solution of $\nameM$ following the procedure described before, i.e., we update $({\tilde{w}^+},{\tilde{w}^-}, \tilde{\bb}, \tilde{\xi}, \tilde{z},\tilde{v})$.
This feasible solution could be improved using the information given by $\tilde{v}$ and $\tilde{z}$ values to obtain the  model below.  Note that constraints associated with a $\tilde{z}$-value equal to 1 are not considered:
\begin{alignat}{4}
(\overline{\mbox{SVM-$\ell_1$}})_{\tilde{v},\tilde{z}} & \quad &  & \mbox{min}&  \quad & 	\displaystyle \sum_{k\in D: \tilde{v}_k=1}\left(w_k^+ + w_k^-\right) + C  \left( \sum_{i\in N:\tilde{z}_i=0}\xi_i \right), & & \nonumber \\
& & &\mbox{s.t.} &\quad &  y_i \left(\sum_{ k \in D: \tilde{v}_k=1}(w_k^+-w_k^-) x_{ik}+\bb \right) \geq 1- \xi_i,     & \quad &i \in N :\tilde{z}_i=0, \nonumber \\
	&&& & &  w_k^+\geq 0, w_k^- \geq 0, &  \quad & k \in D:\tilde{v}_k=1,\nonumber\\ 
& & & & & 0 \leq   \xi_i \leq 2, & \quad  & i \in N:\tilde{z}_i=0 . \nonumber  
\end{alignat}

A feasible solution of \nameM, $({\bar{w}^+},{\bar{w}^-}, \bar{\bb}, \bar{\xi},\bar{v},\bar{z})$, can be obtained considering
the solution of the above linear problem $({\bar{w}^+},{\bar{w}^-}, \bar{\bb}, \bar{\xi})$ together with $\bar{z}=\tilde{z}$ and $\bar{v}=\tilde{v}$.  
From this feasible solution, we compute an upper bound of the model, named UB. Then, we can establish initial bounds of big M parameters as follows:
\begin{eqnarray*}
M_i&=&\underset{j \in N}{\max } \; \left\{ \norm{x_i-x_j}_{\infty} : y_i=y_j \right\}\cdot \mbox{UB},\,\,\mbox{for }i\in N,\\
u_k&=&\mbox{UB},\,\,\mbox{for }k\in D,\\
l_k&=&\mbox{UB},\,\,\mbox{for }k\in D.
\end{eqnarray*}
\subsection{Improving big M parameters for \mbox{\nameM} }

In order to tighten the big M parameters for \nameM, some strategies have been developed. In fact, these strategies are based on solving linear models derived from the original model. We will now detail how each of these parameters could be improved.

\subsubsection{Tightening bounds on $w_k^+$ and $w_k^-$-variables}

By solving the next linear model, we are maximizing the value that the sum of $w_k^+$ and $w_k^-$-variable solutions can have in the feasible region of the original model.
\begin{alignat}{4}
\mbox{(UB-$w$)}	& \quad &  & \max&  \quad & 	\displaystyle \sum_{k=1}^{d}\left(w_k^+ + w_k^-\right) , & & \nonumber \\
	& & &\mbox{s.t.} &  \quad &\eqref{rl-l2-ec2},\eqref{MILP1:ec3},\eqref{UB:w}-\eqref{ec:linealizado},
	& \quad & \nonumber \\
	& & & & &  \sum_{k=1}^{d}\left(w_k^+ + w_k^-\right) + C  \left( \sum_{i=1}^n\xi_i + 2 \sum_{i=1}^n z_i \right)\leq \mbox{UB},    & \quad & \label{aux1} \\
	&&& & &  0\leq v_k\leq 1 & \quad &k\in D,    \label{aux2}  \\
	&&& & &  0\leq z_i\leq 1, & \quad &i\in N.   \label{aux3} 
\end{alignat} 

The optimal objective value of the previous model allows us to obtain an upper bound on $\displaystyle\sum_{k=1}^{d}\left(w_k^+ + w_k^-\right)$. We denote it as UB$_w$. Then, the following valid inequalities are obtained and added to the formulation: 
\begin{equation}\label{ec:boundw}
    w_k^++w_k^-\leq \mbox{UB}_w,\,\,\mbox{for }k\in K.
\end{equation}
Observe that this upper bound tightens the big M parameters, by establishing,
\begin{eqnarray*}
M_i&=&\underset{j \in N}{\max } \; \left\{ \norm{x_i-x_j}_{\infty} : y_i=y_j \right\}\cdot \mbox{UB}_w,\,\,\mbox{for }i\in N,\\
u_k&=&\mbox{UB}_w,\,\,\mbox{for }k\in D,\\
l_k&=&\mbox{UB}_w,\,\,\mbox{for }k\in D.
\end{eqnarray*}
\subsubsection{Tightening bounds on $\bb$-variable}
In a similar way, $\bb$-variable can be tightened by solving
\begin{alignat}{4}
\mbox{(UB-$b$)}/\mbox{(LB-$b$)}	& \quad &  & \max/\min&  \quad & \bb, & & \nonumber \\
	& & &\mbox{s.t.} &  \quad &\eqref{rl-l2-ec2},\eqref{MILP1:ec3},\eqref{UB:w}-\eqref{ec:linealizado},\eqref{aux1}-\eqref{ec:boundw}.& \quad & \nonumber 
\end{alignat} 
The solution of this model allows us to obtain a valid lower and upper bound on the $\bb$-variable of the model. We include them in the formulation of the problem as the following constraint: \begin{equation}
\mbox{LB}_{\bb} \leq \bb \leq \mbox{UB}_{\bb}. \label{eq:b}
\end{equation}

\subsubsection{Tightening bounds on $M_i$ parameters for $i\in N$}

This subsection is focused on tightening the bounds of $M_i$ for $i\in N$. Observe that for $z_i=1$, $i \in N$, constraint \eqref{ec:linealizado} is equivalent to,
\begin{equation}
M_i \geq 1- \xi_i -y_i\left( \sum_{k=1}^{d}(w_k^+ - w_k^-) x_{ik}+\bb\right). \label{newobj}
\end{equation}
We introduce the following model in which the right-hand side of \eqref{newobj} is  maximized over the feasible region of the original model.
\begin{alignat}{4}
	&\left(\mbox{UB}_{M_i}\right) \quad &  & \max&  \quad & 1-\xi_i-y_i\left(\sum_{k=1}^d(w_k^+-w_k^-)x_{ik}+b\right), & & \nonumber \\
	& & &\mbox{s.t.} &  \quad &\eqref{rl-l2-ec2},\eqref{MILP1:ec3},\eqref{UB:w}-\eqref{ec:linealizado},\eqref{aux1}-\eqref{eq:b}.& \quad & \nonumber 
\end{alignat} 
For each $i\in N$, the objective value of this model provides a tightened upper bound on $M_i$. However, in datasets that contain a huge number of individuals, this strategy is computationally demanding. For this reason, we propose the following variant. 
\begin{itemize}
    \item[i)] For all $i\in N,$ when $y_{i}=+1,$ we establish $M_{i}$ as the optimal objective value of the following linear problem:
	\begin{alignat}{4}
	&\left(\mbox{UB}_{M_+}\right)& \quad& \max&  \quad & 	 1- \left(\sum_{k=1}^{d}w_k^+ \underline{x}_{+k} -\sum_{k=1}^{d}w_k^- \bar{x}_{+k}+\bb \right),   \nonumber \\
	&&&\mbox{s.t.} &  \quad &  \eqref{rl-l2-ec2},\eqref{MILP1:ec3},\eqref{UB:w}-\eqref{ec:linealizado},\eqref{aux1}-\eqref{eq:b}, && \nonumber 
	\end{alignat}
	where  $\underline{x}_{+k}=\underset{i \in N.}{\min } \; \{ x_{ik}: y_i=+1 \}$ and $\bar{x}_{+k}=\underset{i \in N.}{\max } \; \{ x_{ik}: y_i=+1 \},$ for $k\in D$.
	
	\item[ii)] For all $i\in N,$ when $y_{i}=-1,$ we establish $M_{i}$ as the optimal objective value of the following linear problem:
	\begin{alignat}{4}
	&\left(\mbox{UB}_{M_-}\right)& \quad& \max&  \quad & 	 1+ \left(\sum_{k=1}^{d}w_k^+ \bar{x}_{-k} -\sum_{k=1}^{d}w_k^- \underline{x}_{-k}+\bb \right),   \nonumber \\
	&&&\mbox{s.t.} &  \quad &  \eqref{rl-l2-ec2},\eqref{MILP1:ec3},\eqref{UB:w}-\eqref{ec:linealizado},\eqref{aux1}-\eqref{eq:b}, && \nonumber
	\end{alignat}
	where  $\bar{x}_{-k}=\underset{i \in N.}{\max } \; \{ x_{ik}: y_i=-1 \}$ and $\underline{x}_{-k}=\underset{i \in N}{\min } \; \{ x_{ik}: y_i=-1 \}$, for $k\in D$. 
	\end{itemize}

The strategy to initialize and tighten the bounds on the big M parameters of the model is summarized in Algorithm \ref{a:strategy} which is described as a pseudocode. 

\begin{algorithm}[htbp]\label{a:strategy}
\setstretch{1.15}
	\KwData{Training sample composed by a set of $n$ individuals with $d$ features.}
	\KwResult{Updated values of $M_i, u_k,$ and $l_k,$ for $i\in N$ and $k\in D$.}

	Solve the problem (SVM-$\ell_1$). From its optimal solution, build and solve ($\overline{\mbox{SVM-}\ell_1}$)$_{\tilde{v},\tilde{z}}.$ From its solution build a feasible solution of $\nameM$ and obtain an upper bound (UB).
		
	\For{$i\in N$}{ 
	$M_i=\underset{j \in N}{\max } \;\{ \norm{x_i-x_j}_ \infty : y_i=y_j \} \cdot \mbox{UB}$.}  
		\For{$k\in D$}{ 
	$u_k=\mbox{UB}, \;\;\; l_k=\mbox{UB}.$} 
		
				Solve \mbox{(UB-$w$)}.
			Let $\mbox{UB}_{w}$ be the optimal objective value of this problem. 
			
			Update $M_i= \underset{j \in N}{\max } \;\{ \norm{x_i-x_j}_ \infty : y_i=y_j \}\cdot \mbox{UB}_{w}, u_k=\mbox{UB}_{w},$ and $l_k=\mbox{UB}_{w}.$ Add the obtained bounds to the formulation $\nameM$ including the set of constraints  \eqref{ec:boundw}.\\
			Obtain $\mbox{LB}_{\bb}$ and $\mbox{UB}_{\bb}$, lower and upper bounds respectively, of the $\bb$-variable by solving \mbox{(LB-$b$)} and \mbox{(UB-$b$)}.\\
Include the constraint \eqref{eq:b} in the formulation $\nameM$.

\Do{an improvement of the bounds is obtained}{

Repeat Steps 6 and 7 including constraints \eqref{ec:boundw} and \eqref{eq:b} in model $(\mbox{UB-}w)$.

		\Case{\text{Variant I}}{
					\For{$i\in N$}{ Update $M_{i}$ as the optimal value of the problem  $\left(\mbox{UB}_{M_{i}}\right)$.
			}
		}
		\Case{\text{Variant II}}{ 
				For $i\in N,$ when $y_i=1,$ update $M_i$ as the optimal value of the problem $\left(\mbox{UB}_{M_+}\right)$.\\
			For $i\in N,$ when $y_i=-1,$ update $M_i$  as the optimal value of the problem $\left(\mbox{UB}_{M_-}\right)$.
		}
}
			\caption{Variant 1 and 2. Computation of big M parameters.}
\end{algorithm}

To summarize, in this section we have introduced a model based on SVM using the $\ell_1$-norm that includes feature selection and avoids the effect of outliers. For this model, we have developed a mixed integer linear formulation $\nameM$. In addition, we have introduced a methodology for initializing and tightening the big M parameters appearing in the model. In the next section, we present a heuristic approach to obtain adequate feasible solutions for the model.

\section{Heuristic approach}
\label{sec:Heuristic}
In the previous section, some strategies to tighten big M parameters of $\nameM$ have been presented. Consequently, the formulation is strengthened by the obtained bounds easing its exact solution approach. However, the exact solution of this model cannot be obtained in reasonable times for large datasets. For this reason, it is necessary to develop a heuristic solution method to deal with such datasets.

In this section, we develop a procedure for obtaining good quality solutions. It is based on the Adaptative Kernel Search (AKS) proposed by \cite{AKS}. The idea of the AKS is to solve a sequence of mixed integer problems derived from the original one, where the majority of the binary variables are fixed and only some of them  are considered as integer -- this set of integer variables is called Kernel.

The heuristic approach that we propose for $\nameM$ has two main characteristics. On the one hand, we introduce an AKS for $v$-variables, see \citet{22Luisa}. On the other hand, the $z$-variables are dynamically fixed using the solution of the previous sub-problems. This strategy considers sub-problems where the number of binary variables and constraints containing the big M parameter are reduced with respect to the
original problem. As a consequence, these sub-problems have a shorter solution time. 
This heuristic method will be referred to as Dynamic Adaptative Kernel Search (DAKS). DAKS allows us to obtain adequate feasible solutions of $\nameM$ in datasets of large size, as can be appreciated in Section \ref{sec:computationalExp}. The improvement is justified by the reduction in the number of integer variables and big M constraints in the sub-problems considered in this procedure.

\subsection{Dynamic Adaptative Kernel Search (DAKS)}

The DAKS algorithm has three phases. 
\begin{itemize}
    \item The first one computes an initial feasible solution to the problem. Using this solution, the big M parameters appearing in the model are initialized and some $z$-variables are fixed. 
    \item The second one determines the initial set of integer $v$-variables; i.e., the initial kernel. Furthermore,  in this step, the remaining $v$-variables are ranked. Note that the most promising  variables to take a value of 1 in the optimal solution are considered first. The chosen order has a huge influence on the solution obtained by the heuristic, so this order should be adapted to the dataset under study. Observe that although in the original algorithm several sets of integer variables can be considered, in this version, we only apply the Adaptative Kernel Search over $v$-variables, i.e., $z$-variables are adjusted using a different methodology during the procedure.  

\item In the third phase, a sequence of restricted MILP derived from the original problem are solved, including constraints that ensure that a progressively better bound on the solution is obtained. The solution and the computational effort required to solve the previously restricted MILP, guide the construction of the subsequent Kernels. In each iteration, the
Kernel is updated to test other  promising variables and to remove useless ones. Moreover, taking into account the solution of the previous sub-problem, each $z$-variable is fixed to zero, fixed to one or considered as a binary variable, i.e., this set of variables is updated dynamically. After some iterations, the second phase is repeated with the objective of reordering $v$-variables using the information of the current $z$-variables values.
\end{itemize}

Given $\mathcal{K}\subseteq D,$ the restricted MILP derived from the original problem will be referred to as $\mbox{\nameM($\mathcal{K}$)}_{\hat{z}},$ where if $k\in \mathcal{K},$ $v_k$-variable is considered a binary variable, it is otherwise fixed to zero. Furthermore, with the objective of simplifying the notation of these sub-problems, we include an auxiliary vector $\hat{z}$ which indicates the corresponding  $z$-variables values. Each element of vector $\hat{z}$ can take one of the following values with the following meanings:
$$
\hat{z}_i=
\begin{cases}
0,&\mbox{it indicates that $z_i$ is fixed to $0$ in the sub-problem,}\\
1,&\mbox{it indicates that $z_i$ is fixed to $1$ in the sub-problem,}\\
2,&\mbox{it indicates that $z_i$ is a binary variable in the sub-problem.}\\
\end{cases}
$$
Therefore, $\mbox{\nameM($\mathcal{K}$)}_{\hat{z}}$ model  is formulated as described below.
\begin{alignat}{4}
& \quad &  & \mbox{min}&  \quad & 	\displaystyle \sum_{k\in \mathcal{K}}\left(w_k^+ + w_k^-\right) + C  \left( \sum_{i\in N:\hat{z}_i\neq 1}\xi_i + 2 \sum_{i\in N:\hat{z}_i=2} z_i +2 \sum_{i\in N:\hat{z}_i=1} 1\right), & &\nonumber \\
& & &\mbox{s.t.} &  \quad & y_i\left( \sum_{k\in \mathcal{K}}\left(w_k^+ - w_k^-\right)x_{ik}+\bb\right) \geq 1- \xi_i -M_iz_i,   & \quad  & i\in N:\hat{z}_i=2, \nonumber\\
& & & &  \quad & y_i\left( \sum_{k\in \mathcal{K}}\left(w_k^+ - w_k^-\right)x_{ik}+\bb\right) \geq 1- \xi_i,   & \quad  & i\in N:\hat{z}_i=0, \nonumber\\
& & & &  \quad &	\xi_i \leq 2(1 - z_i), && i \in N:\hat{z}_i=2, \nonumber\\
		& & &&& \sum_{k\in \mathcal{K}}v_k \leq B, & \quad   & \nonumber  \\
			&&& & & 0 \leq  w_k^+ \leq u_k v_k, & \quad &k\in \mathcal{K},    \nonumber \\
	&&& & &  0 \leq w_k^- \leq l_k v_k, & \quad &k\in \mathcal{K},    \nonumber\\
& & &&& v_k \in \{0,1\}, &\quad   & k\in \mathcal{K},  \nonumber\\
& & & &  \quad & 0 \leq   \xi_i \leq 2, & \quad  & i \in N:\hat{z}_i\neq 1,  \nonumber\\
& & & & & z_i\in \{0,1\}, & &  i \in N:\hat{z}_i=2.  \nonumber
\end{alignat}

 In the next subsections, the three main phases of DAKS procedure are described in detail and a pseudocode of this heuristic approach is also presented in Algorithm \ref{a:ks}.

\subsubsection{Initial phase}
In the first phase,  Algorithm \ref{a:strategy} is applied with the aim of obtaining an initial upper bound of the problem and valid values for the big M parameters. We can set a time limit for the algorithm or establish a fixed number of iterations. The modeler should decide the  best strategy to follow for the dataset under study and choose one of the variants of this algorithm, taking into account the number of individuals and features  considered in the dataset.

Furthermore, in this phase, $z$-variables are fixed. Let $({\bar{w}^+},{\bar{w}^-}, \bar{\bb}, \bar{\xi}, \bar{v},\bar{z})$ be the initial solution built in Algorithm \ref{a:strategy}. Then, vector $\hat{z}$ is determined as follows, 
$$\hat{z_i}=\begin{cases}
0,&\mbox{if $\bar{z}_i=0$ and $\bar{\xi}_i<1$,}\\
1,&\mbox{if $\bar{z}_i=1$},\\
2,&\mbox{if $\bar{z}_i=0$ and $\bar{\xi}_i\geq1$,}\\
\end{cases}$$
for $i\in N$. Note that if $\hat{z}_i$ is established as two, then $z_i$ will be a binary variable in  $\mbox{\nameM($\mathcal{K}$)}_{\hat{z}}$.

The idea of this procedure is to discard as outliers the well-classified individuals in the initial solution, i.e., $\hat{z}_i$ is fixed to zero. Furthermore, the individuals 
whose corresponding $\bar{z_i}$-values are one in the initial solution are set as outliers, i.e., $\hat{z}_i$ is fixed to one. Finally, we establish as ``possible outliers" (considering the  associated $z$-variable as binary) 
the wrong-classified individuals
  corresponding  to $\bar{\xi}$-values in the interval [1,2].

\subsubsection{Second phase}\label{sec:secpha}

In the second phase, we will create an initial kernel on $v$-variables and the remaining variables will be ordered according to how likely these variables
will take a value of one in the optimal solution. This stage is essential for the success of the algorithm. 

We propose the following strategy for ordering the variables. First, we solve the problem $\mbox{\nameM($D$)}_{\hat{z}}$ considering $v$-variables as continuous. We will denote the previously described model as  $\mbox{(Rel$_v$-RL-FS-M)($ D$)}_{\hat{z}}$. Observe that in this relaxed problem, only the $z$-variables whose corresponding $\hat{z}$ take value 2 are binary (recall that $\hat{z}$ was computed in the initial phase). The optimal objective value of this problem will be called $\mbox{LB}_{\hat{z}}$. Again, we solve it by fixing the values of the binary variables in order to obtain the reduced costs. 

 Let $\breve{w}^+_k$ and $\breve{w}^{-}_k$ be
 the optimal values  of variables $w^+_k$ and $w^-_k$  and $r_k^+$ and $r_k^-$ the corresponding reduced cost of these variables.
Hence, the variables are ordered in non-decreasing order with respect to vector $r$, which is computed as:
$$r_k=\begin{cases} 
    -(\breve{w}_k^++\breve{w}_k^-), &\text{if }\breve{w}_k^++\breve{w}_k^->0, \\
\min\{r_k^+,r_k^-\},   &\text{otherwise}
\end{cases}$$

 The initial kernel, $\mathcal{K}_0$, is composed by $k\in D$ such that $w^+_k$ or $w^-_k$ variable takes a positive value in the solution of  $\mbox{(Rel$_v$-RL-FS-M)($D$)}_{\hat{z}}$ model considering $z$-variables fixed to their optimal values in $\mbox{(Rel$_v$-RL-FS-M)($D$)}_{\hat{z}}$.
 Moreover, the first time that this phase is executed, $\mathcal{K}_0$ includes the indexes $k\in D$ such that $w^+_k$ or $w^-_k$ variables take a positive value in the initial solution built in the first step of Algorithm \ref{a:strategy}.

 Once this initial kernel $\mathcal{K}_0$ is defined, the model $\mbox{\nameM($\mathcal{K}_0$)}_{\hat{z}}$ is solved. Its solution is an upper bound (UB) on the problem \nameM. Let $({\bar{w}^+},{\bar{w}^-}, \bar{\bb}, \bar{\xi}, \bar{v},\bar{z})$ be the solution of the problem, i.e., the solution of the current upper bound. After this step, vector $\hat{z}$ is updated as follows: 
 \begin{itemize}
     \item If $\hat{z}_i$ is equal to zero and $\bar{\xi}_i>1$, for $i\in N$, (i.e., the individual was discarded as outlier, but in the current solution it is wrongly classified), $\hat{z}_i$ will be fixed to two (i.e., it will be set as a ``possible outlier'').
     \item If $\hat{z}_i$ is equal to one and $y_i\left( \sum_{k\in D}\left(\bar{w}_k^+ - \bar{w}_k^-\right)x_{ik}+\bar{\bb}\right) \geq 0$, for $i\in N$, (i.e., the individual was set as an outlier, but in the current solution it is well classified), $\hat{z}_i$ is updated to two (i.e., it will be set as a ``possible outlier'').
 \end{itemize} 

Note that the criterion of being a ``possible outlier'' can be adjusted to the dataset, imposing tighter conditions, e.g. $\bar{\xi}_i>0$, or relaxing them, e.g. $\bar{\xi}_i>1.5$.

\subsubsection{Third phase}
An iterative procedure starts in this phase. In each iteration, ($it$), a new set of indexes, named $B_{it}\subseteq D$ is added to the kernel. These indexes are included in the order determined by vector $r$, i.e., the most promising $v$-variables to take value one in the optimal solution are considered first. Note that each index is included at most once. The initial size of $B_{it}$ is the size of $\mathcal{K}_0$, see \cite{Guastaroba2012,22Luisa}. In the subsequent iterations, the size of $B_{it}$ will be adjusted.

The model $\mbox{\nameM($\mathcal{K} \cup B_{it}$)}_{\hat{z}}$ is solved adding the following constraints: 
\begin{eqnarray}
&&\sum_{k\in \mathcal{K} \cup B_{it}}(w_k^+ + w_k^-) + C  \left( \sum_{i=1}^n\xi_i + 2 \sum_{i=1}^n z_i \right)\leq \mbox{UB},  \label{betterUB}\\
&&\sum_{k\in B_{it}}v_k \geq 1. \label{atleast1}
\end{eqnarray}

Constraint \eqref{betterUB} ensures that a better upper bound of the problem is obtained, and constraint \eqref{atleast1} ensures that at least one variable is selected for the new $B_{it}$ set. If the previous iteration found a feasible solution but the optimality of this solution was not proved because the time limit was reached, constraint \eqref{atleast1} will be replaced  by: 
\begin{equation}
\sum_{k\in B_{it}  \cup \{k \in \mathcal{K}: \bar{v}_k=0\}}v_k \geq 1, \label{atleast1-no}
\end{equation}
where $\bar{v}_k$ represents the solution associated with the current upper bound.

A feasible solution of this sub-problem, i.e., $\mbox{\nameM($\mathcal{K}\cup B_{it}$)}_{\hat{z}}+\eqref{betterUB}+\eqref{atleast1}$ if the previous solution was optimal or $\mbox{\nameM($\mathcal{K}\cup B_{it}$)}_{\hat{z}}+\eqref{betterUB}+\eqref{atleast1-no}$ if the previous solution was feasible, will improve the upper bound. However, this problem may not be feasible and for this reason, a time limit $t_{limit}$ is imposed. Furthermore, we also add a time limit $t_{fea}$ for finding a feasible solution, i.e, the sub-problem is stopped if no feasible solution is found within $t_{fea}$. In addition, as the big M parameters affect the lower bounds of the problem, proving optimality is hard. Consequently, if the incumbent solution is not improved after a fixed time $t_{inc}$, the resolution of the sub-problem is stopped.

After this, the kernel is updated: the selected variables of the set $B_{it}$ are added, denoted as $\mathcal{K}^{+}$, and the variables of the kernel that were not selected in the last $p$ iterations are removed from the kernel, denoted as $\mathcal{K}^-.$ 
Moreover, vector $\hat{z}$ is updated as follows: 
\begin{table}[htb]
\begin{center}
\begin{tabular}{|c|c|c|}
 \hline
                  Previous value of vector $\hat{z}$&  Previous and current solutions & Updated value of vector $\hat{z}$  \\
                   \hline
                  $\hat{z}_i$=0& Current solution $\bar{\xi}_i\geq1$ &   $\hat{z}_i$=2\\
                   \hline
            $\hat{z}_i$=1 & $y_i\left( \sum_{k\in D}\left(\bar{w}_k^+ - \bar{w}_k^-\right)x_{ik}+\bar{\bb}\right) \geq 0$ & $\hat{z}_i$=2  \\
                   \hline
\multirow{2}{*}{$\hat{z}_i$=2} & Last $q$ solutions $\bar{z}_i=1$ &  $\hat{z}_i$=1 \\
\cline{2-3}
                  & Last $q$ solutions $\bar{z}_i=0$ &  $\hat{z}_i$=0 \\
                   \hline
\end{tabular}
\end{center}
\caption{Procedure to update vector $\hat{z}$.}
\label{tab1}
\end{table}

 The idea of the previous procedure is: if the individual was discarded as outlier  but in the current solution it is wrongly classified, it will be set as a ``possible outlier''. If the individual was set as outlier, but in the current solution it is well classified, it will be set as a ``possible outlier''. Finally, if the individual was set as ``possible outlier''  but it has the same value in the last $q$ solutions, the variable will be fixed to this value. Observe that the criterion of being a ``possible outlier'' can be adjusted to the dataset, imposing tighter conditions or relaxing them, as explained at the end of Section \ref{sec:secpha}.

In the following iterations, the problem  $\mbox{\nameM($\mathcal{K} \cup B_{it}$)}_{\hat{z}}$ is solved. The time employed in solving the previous iteration determines the size of  $\mathcal{K} \cup B_{it}$, $S_{K \cup B_{it}}$. If the problem related to an iteration is easily solved, i.e., if $t_{K \cup B_{it}} \leq t_{Easy}$, the number of $v$-variables that are binary in the next step is incremented, i.e., $S_{K \cup B_{it+1}}=(1+\delta)S_{K \cup B_{it}}$, where $0\leq \delta \leq 1.$

This iterative process (phase 3) will be stopped if any of the following situations occur: 
\begin{enumerate}[i)]
    \item If $\mbox{LB}_{\hat{z}}=\mbox{UB}$ and $\hat{z}$ remains unchanged, the algorithm finishes. If $\mbox{LB}_{\hat{z}}=\mbox{UB}$, we know that the obtained solution cannot be improved having fixed $z$-variables to these values. Since different $\hat{z}$-values to the current ones are not proposed for the next iteration, the algorithm finishes. Note that $\mbox{LB}_{\hat{z}}$ is an upper bound of the lower bound (LB) of the problem. Hence, the solution could be a local minimum or the optimal solution of the problem.
     \item If $\mbox{LB}_{\hat{z}}=\mbox{UB},$  the current $\hat{z}$ is the same as the one used to compute $\mbox{LB}_{\hat{z}}$, and different values are proposed for the next iteration to update $\hat{z}$ in the procedure described in Table \ref{tab1}; return to phase two. If this situation happens, our solution cannot be improved if $z$-variables are fixed to these values. However, it may improve if $z$-variables are fixed to different values. In order to recompute the initial kernel and  sort the $v$-variables considering the updated values of $\hat{z}$, the algorithm returns to phase two.
     \item After a criterion defined by the user (e.g. doing a fixed number of iterations, establishing a time limit or fixing a percentage of $v$-variables that should be analyzed), the algorithm  returns to phase two. The objective is to sort the $v$-variable taking into account the current $\hat{z}$.
\end{enumerate}

Note that the number of iterations that the algorithm does in Phase 2 and 3 is a criterion determined by the user that should be adapted to the dataset. Similarly, the proposed heuristic has many parameters $p,q,\delta,t_{limit},t_{Easy}, \dots$ that the user should adapt to the dataset and their needs, finding the desired balance between time employed and required precision.
 In Algorithm~\ref{a:ks}, this procedure is summarized and described as a pseudocode. 

In conclusion, in this section we have proposed a new algorithm to obtain accurate feasible solutions of our model based on the kernel search algorithm. In particular, we have developed a different strategy to  seek the most promising values of $z$-variables with the objective of avoiding using the reduced costs of these variables  due to the high influence that the big M parameters have over them.

\begin{algorithm}[htbp]\label{a:ks}
	\KwData{Training sample composed by a set of $n$ individuals with $d$ features.}
	\KwResult{Heuristic solution of \nameM.}

	Initialize big M parameters applying Algorithm \ref{a:strategy}. Let $({\bar{w}^+},{\bar{w}^-}, \bar{\bb}, \bar{\xi}, \bar{v},\bar{z})$ be the initial solution built in this Algorithm.
	
	\For{$i\in N$}{ 
\If{$\bar{z}_i=1$}{$\hat{z}_i$=1}
\If{$\bar{z}_i=0$ and $\bar{\xi}_i>1$}{$\hat{z}_i$=2}
\If{$\bar{z}_i=0$ and $\bar{\xi}_i\leq1$}{$\hat{z}_i$=0}
}

\Do{a criterion is fulfilled (e.g. number of iterations, time limit)}{

Solve $\mbox{(Rel$_v$-RL-FS-M)($D$)}_{\hat{z}}$. Solve it again fixing the values of the binary $z$-variables. Sort $v$-variables in non-decreasing order with respect to vector $r$. 

Constitute  initial $\mathcal{K}=\mathcal{K}_0$. Solve $\mbox{\nameM($\mathcal{K}$)}_{\hat{z}}$. Let $({\bar{w}^+},{\bar{w}^-}, \bar{\bb}, \bar{\xi}, \bar{v},\bar{z})$ be the solution of the current UB. 

	\For{$i\in N$}{ 
\If{($\hat{z}_i=0$  \textbf{and} $\bar{\xi}_i>1$) \textbf{or} ($\hat{z}_i$=1 \textbf{and} $y_i\left( \sum_{k\in D}\left(\bar{w}_k^+ - \bar{w}_k^-\right)x_{ik}+\bar{\bb}\right) \geq 0$)}{$\hat{z}_i$=2}
}  
$it=0$

\Do{a criterion is fulfilled (e.g. number of iterations, time limit)}{
\If{$it>0$ and the solution of $it-1$ was feasible but not optimal}
{

 Let $B_{it}$ be the new set of promising $v$-variables related to iteration $it$.

Solve $\mbox{\nameM($\mathcal{K}\cup B_{it}$)}_{\hat{z}}+\eqref{betterUB}+\eqref{atleast1-no}$, where $B_{it}$ is the set of promising $v$-variables associated with iteration $it$.
}
\Else{
Solve $\mbox{\nameM($\mathcal{K}\cup B_{it}$)}_{\hat{z}}+\eqref{betterUB}+\eqref{atleast1}$.
}

\If{an optimal or feasible solution is obtained}
{Update $({\bar{w}^+},{\bar{w}^-}, \bar{\bb}, \bar{\xi}, \bar{v},\bar{z})$ and UB.

Let:
\begin{itemize}
\item[]$\mathcal{K}^{+}:=\{k\in B_{it}: \bar{v}_k=1 \}.$ 
\item[]$\mathcal{K}^{-}:=\{k\in \mathcal{K}: k \mbox{ has not been selected in the last $p$ feasible iterations.} \}.$
\end{itemize}
Update $\mathcal{K}=\mathcal{K}\cup\mathcal{K}^{+}\setminus\mathcal{K}^-$ and $\hat{z}$ as described in Table \ref{tab1}.}

\If{$\mbox{LB}_{\hat{z}}=\mbox{UB}$ and $\hat{z}$ remains unchanged} {Go to Line 31.}
 \If{$\mbox{LB}_{\hat{z}}=\mbox{UB}$ and $\hat{z}$ is updated} {Go to Line 9.}
 
 $it=it+1$.

}
}
\Return $({\bar{w}^+},{\bar{w}^-}, \bar{\bb}, \bar{\xi}, \bar{v},\bar{z}).$
	\caption{Dynamic Adaptative Kernel Search (DAKS).}
\end{algorithm}
\section{Computational Experiments}
\label{sec:computationalExp}

In this section, we set out the results of several computational experiments. Firstly, we analyze the performance of  Algorithm~\ref{a:ks} comparing the solution provided using this algorithm to the solution obtained applying an exact resolution method. The results demonstrate that the proposed algorithm computes high-quality solutions. Secondly, we present a comparison of the classification done by several models over real-life datasets showing the efficiency of the proposed classifier. 
  
The experiments were conducted on an Intel(R) Xeon(R) W-2135 CPU 3.70 GHz 32 GB RAM, using CPLEX 12.9.0. in  Concert Technology C++. As seen in \cite{bonami} due to the presence of big M parameters, the relative MIP tolerance and the integrality tolerance were fixed to zero. The remaining parameters were left to their default values unless stated otherwise.

\subsection{Data}
The computational experiments were carried out on real-life datasets.  They are specified in Table \ref{tab:RWDS}, where $n$ is the number of individuals, $d$ is the number of features, and the last column reports the percentage of elements in each class. Leukemia dataset is from \cite{Leukemia}, Colon dataset is from \cite{colon}, and DLBCL dataset is from \cite{DLBCL}. All other datasets are from the UCI repository \citep[see][]{DatasetUCI}. Observe that these datasets were previously used to analyze the performance of models based on SVM \citep[see, for instance,][]{31Brooks2011,22Luisa,21Maldonado2014}.

\begin{table}[htbp]
	\resizebox{\hsize}{!}{
	\centering
	\begin{tabular}{llrrr}
		\toprule
		Label & Name in repository  & $n$     & $d$     & Class(\%) \\
		\midrule
		Colon & Analysis of tumor and normal colon tissues& 62 &2000& 35/65 \\
		Leukemia   &  Gene expression measurements on leukemia patients & 72   & 5327    &  47/53\\
		DLBCL   &  Measurements from biopsies of Diffuse Large B-cell Lymphoma patients& 77   & 7129    &  75/25\\
		
		SONAR  & Connectionist bench (sonar, mines vs. rocks) & 208   & 60     & 54/46 \\
		IONO & Ionosphere & 351   & 33    & 64/36 \\
		
    	Arrhythmia   & Cardiac arrhythmia & 420   & 258    & 57/43 \\
		Wdbc   & Breast cancer Wisconsin (diagnostic) & 569   & 30    & 63/37 \\
		Mfeat &  Multiple Features -  Handwritten numerals & 2000&649& 10/90\\
		Lepiota& 
Mushroom - Species of mushrooms in the Agaricus and Lepiota Family& 8124& 109& 52/48\\
		\bottomrule
	\end{tabular}%
	}
	\caption{Real-life datasets}
	\label{tab:RWDS}%
	
\end{table}%

In these datasets, the following values for the parameter $C$ will be analyzed, $C \in \{0.01,0.1,1,$ $10,100\},$ as in \cite{31Brooks2011}. We tested five different values for the parameter $B$, which is different for each dataset. To choose it, we solved \mbox{\nameM} on ten randomly generated samples containing 90\% of individuals from the original dataset being $B$ equal to the number of features. Note that in each sample, the percentage of elements from each class is the same as in the original dataset. We compute the average from the number of selected features for each $C$. The maximum  average is the greatest $B$ analyzed, the other four are proportional to it, being $\frac{2}{3},\frac{1}{2},\frac{1}{3},$ and $\frac{1}{5}$ of the greatest value of $B$. In Table \ref{tab:Bvalues}, these values are depicted. Moreover, since 
the two greatest values of $B$ for the Arrhythmia dataset ($161$ and $107$) and Mfeat dataset ($404$ and $269$) computed by the above description are too big, we have removed them and we have included the following values of $B$: $16$ and $8$ for Arrhythmia dataset; $30$ and $10$ for Mfeat dataset. 
\begin{table}[htbp]

	\centering
	\begin{tabular}{ll}
		\toprule
		Label & $B$  \\
		\midrule
		Colon & $\{34,23,17,11,7\}$ \\
		Leukemia   & $\{48,32,24,16,10\}$\\
		DLBCL   &  $\{44,29,22,15,9\}$\\
		SONAR  & $\{57,38,29,19,11\}$\\
		IONO & $\{33,22,17,11,7\}$ \\
		
    	Arrhythmia   & $\{81,54,32,16,8\}$\\
		Wdbc   & $\{29,19,15,10,6\}$ \\
		Mfeat &  $\{202,135,81,30,10\}$\\
		Lepiota& $\{33,22,17,11,7\}$\\
		\bottomrule
	\end{tabular}%
	\caption{Values of parameter $B$.}
	\label{tab:Bvalues}%
	
\end{table}%

\subsection{Validation of DAKS algorithm}

In this section, we validate the  efficiency of Algorithm \ref{a:ks}, i.e, we compare the solution obtained by the exact method with the one provided by Algorithm \ref{a:ks}. The exact method consists of applying Algorithm \ref{a:strategy} and solving the resulting formulation using CPLEX. The resolution of this model is quite hard, and for this reason, a time limit of 7200 seconds has been set. 

We solved $\nameM$  for all datasets for all $C$ values ($C \in \{0.01,0.1,1,10,100\}$) and the smallest value of $B$ tested, i.e., Colon $B=7$, Leukemia $B=10$, DLBCL $B=9$, SONAR $B=11$, Iono $B=7$, Arrhythmia $B=8$, Wdbc $B=6$, Mfeat $B=10$, and Lepiota $B=7$.  For the Mfeat and Lepiota datasets, we applied Variant 2 of Algorithm \ref{a:strategy} to initialize the big M parameters, while Variant 1 was used in the rest of the cases. 

As done in \cite{AKS} for the AKS algorithm, parameter $\delta$ was set to 0.35 in DAKS. Moreover, we established $p=2,$ $q=2,$ $t_{Easy}=10,$ $t_{fea}=120,$ and $t_{inc}=160$. Furthermore, in the restricted problems, we limited the computational time ($t_{limit}$) to 400 seconds. In order to improve the behavior, as in \cite{AKS}, we used CPLEX with default values, except for the following parameters: BndStrenInd=1, MIPEmphasis=HiddenFeas, FPHeur=1. The parameters previously mentioned to avoid the negative effect of big M values were also used. 

The results are depicted in Table \ref{tab:validation}. The first column indicates the name of the dataset and the second column the value of parameter $C$. The next column reports the total time of the exact solution method in seconds, i.e., the strategy time plus the solving time. The following column contains information about the MIP relative GAP reported by CPLEX. The fifth column indicates the time in seconds that the heuristic took. Finally, the sixth column reports the percentage difference between the best solution found by heuristic algorithm ($\text{BS}_{h}$) the best solution found by the exact solution method ($\text{BS}_{e}$). This difference is computed as follows:  $$\%\text{ BS}=\dfrac{\text{ BS}_{h}-\text{ BS}_{e}}{\text{ BS}_{e}}.$$ 

\begin{table}[htb]
  \centering
  \resizebox{\hsize}{!}{
    \begin{tabular}{c|r|rrrr}
    \multicolumn{1}{l|}{Data} & \multicolumn{1}{r|}{C} & \multicolumn{1}{r}{$t_{e}$} & \multicolumn{1}{r}{GAP} & \multicolumn{1}{r}{$t_{h}$} & \multicolumn{1}{r}{\%BS} \\
    \hline
    \multirow{5}[2]{*}{\begin{turn}{90}Colon \end{turn}} & 100   & 7225.01 & 23.59\% & 49.06 & 0.00\%  \\
          & 10    & 7219.15 & 21.95\% & 34.66 & 2.40\% \\
          & 1     & 7218.79 & 5.35\% & 15.26 & 1.54\%   \\
          & 0.1   & 2856.23 & 0.00\% & 6.69  & 0.00\%  \\
          & 0.01  & 1.25  & 0.00\% & 1.24  & 0.00\%  \\
    \hline
    \multirow{5}[2]{*}{\begin{turn}{90}Leukemia\end{turn}} & 100   & 7230.38 & 11.45\% & 120.14 & 3.23\%  \\
          & 10    & 7241.98 & 12.77\% & 101.67 & 2.05\%  \\
          & 1     & 7281.67 & 14.82\% & 92.79 & -1.90\%  \\
          & 0.1   & 41.10 & 0.00\% & 43.90 & 0.00\%  \\
          & 0.01  & 4.75  & 0.00\% & 4.71  & 0.00\%  \\
    \hline
    \multirow{5}[2]{*}{\begin{turn}{90}DLBCL\end{turn}} & 100   & 7250.71 & 10.42\% & 73.15 & 0.37\% \\
          & 10    & 7266.44 & 10.76\% & 82.42 & 0.37\% \\
          & 1     & 7326.72 & 10.35\% & 113.93 & 2.33\%  \\
          & 0.1   & 2437.37 & 0.00\% & 58.74 & 0.00\%  \\
          & 0.01  & 7.10  & 0.00\% & 7.36  & 0.00\% \\
    \hline
    \multirow{5}[2]{*}{\begin{turn}{90}SONAR\end{turn}} & 100  
  &    7217.41 & 94.46\% & 302.70 & 4.83\% \\
&  10 &  7211.48 & 90.50\% & 202.88 & 1.91\% \\
& 1 &  7205.42 & 59.18\% & 578.35 & 0.56\% \\
          & 0.1   & 7207.64 & 55.09\% & 1487.38 & 0.00\%  \\
          & 0.01  & 7207.57 & 1.02\% & 2.76  & 0.00\%  \\
    \hline
    \multirow{5}[2]{*}{\begin{turn}{90}IONO\end{turn}} & 100   & 7222.98 & 84.66\% & 77.35 & -8.55\%  \\
          & 10    & 7226.58 & 72.06\% & 173.36 & 0.00\%  \\
          & 1     & 7237.98 & 63.71\% & 117.58 & -3.07\% \\
          & 0.1   & 7221.80 & 46.56\% & 305.99 & -0.07\% \\
          & 0.01  & 7222.49 & 12.13\% & 378.45 & 0.00\% \\
    \hline
     \end{tabular}%
     \hspace{0.4cm}
     \begin{tabular}{c|r|rrrr}
     \multicolumn{1}{l|}{Data} & \multicolumn{1}{r|}{C} & \multicolumn{1}{r}{$t_{e}$} & \multicolumn{1}{r}{GAP} & \multicolumn{1}{r}{$t_{h}$} & \multicolumn{1}{r}{\%BS}  \\
     \hline
    \multirow{5}[2]{*}{\begin{turn}{90}Arrhythmia\end{turn}} &   100   & 7244.87 & 98.89\% & 1844.90 & -10.00\% \\
   & 10    & 7268.66 & 98.38\% & 1860.49 & -8.91\% \\
  &  1     & 7246.76 & 89.60\% & 1844.65 & -1.91\% \\
    &0.1   & 7255.37 & 83.78\% & 1832.81 & 0.19\% \\
    &0.01  & 7258.67 & 59.64\% & 537.19 & 0.00\% \\
    \hline
    \multirow{5}[2]{*}{\begin{turn}{90} Wdbc\end{turn}} & 100   &      1954.27 &     0.00\%  &  148.65     &    1.92\%    \\
          & 10    & 204.05 & 0.00\% & 44.96 & 0.00\%  \\
          & 1     & 41.17 & 0.00\% & 39.31 & 0.00\%  \\
          & 0.1   & 29.74 & 0.00\% & 29.17 & 0.00\%  \\
          & 0.01  & 7221.52 & 16.27\% & 378.81 & 0.00\%  \\
    \hline
    \multirow{5}[2]{*}{\begin{turn}{90} Mfeat\end{turn}} &     100   & 7223.05 & 9.31\% & 131.33 & 2.21\% \\
    &10    & 7218.46 & 8.04\% & 96.09 & 2.21\% \\
    &1     & 7228.61 & 7.75\% & 89.22 & 0.80\% \\
    &0.1   & 1060.93 & 0.00\% & 70.27 & 0.00\% \\
    &0.01  & 7233.85 & 57.84\% & 101.85 & 0.00\% \\
    \hline
    \multirow{5}[2]{*}{\begin{turn}{90} Lepiota\end{turn}} & 100   & 63.53 & 0.00\% & 26.52 & 0.00\%  \\
          & 10    & 41.76 & 0.00\% & 29.15 & 0.00\%  \\
          & 1     & 55.67 & 0.00\% & 30.75 & 0.00\%  \\
          & 0.1   & 7249.02 & 11.64\% & 58.62 & 0.00\%  \\
          & 0.01  & 7403.56 & 1.75\% & 276.12 & 0.00\% \\
    \hline
    \end{tabular}%
    }
    \caption{Comparison between the exact resolution method and the heuristic procedure.}
  \label{tab:validation}%
\end{table}%
As can be appreciated in Table \ref{tab:validation}, the proposed heuristic resolution method reported the same solution as the exact solution method in most cases, with the advantage that the heuristic algorithm is much less time-consuming than the exact resolution method. Note also, that in some cases, Algorithm \ref{a:ks} found a better solution than the exact solution method.
For instance, in the case of Leukemia dataset with $C=1$, Algorithm \ref{a:ks} found a 1.9\% better solution in one minute and a half than the exact solution method found in over two hours. Besides, in the case of IONO dataset with $C=100$, Algorithm \ref{a:ks} found a 8.55\% better solution in less than two minutes than the exact solution method found in over two hours. Similarly, in the case of Arrhythmia dataset with $C=100$ and $C=10$, Algorithm \ref{a:ks} found 10\% and 8.91\% better solutions respectively than the exact solution method found in over two hours.

In summary, we compared the performance of Algorithm \ref{a:ks} with an exact solution approach for a challenging $B$ parameter value. We proved that DAKS provides quite good solutions in the 45 cases analyzed for different $C$ parameter values. On average, the difference between the best solution found by the exact solution method and the heuristic is -0.17\%, while the heuristic has significantly reduced the amount of time required. 

\subsection{Model validation}

The aim of this section is to compare the performance of the proposed classifier $\nameM$  with other well-known methods based on SVM in the literature. We focus our attention on efficient models that deal with outliers detection $(\mbox{RL-$\ell_1$-M})$ or feature selection: $(\mbox{FS-SVM})$, $\mbox{(Fisher-SVM)}$, and (RFE-SVM). More concretely, the purpose of $(\mbox{RL-$\ell_1$-M})$  is to avoid the  influence of outliers, see \cite{31Brooks2011}. We used the methodology proposed in \cite{BalMarRod19} to solve this model, imposing 1800 seconds of time limit. On the other hand, the objective of $(\mbox{FS-SVM})$ is to limit the number of features selected by the classifier and it was solved by applying the heuristic algorithm proposed in \cite{22Luisa}. Similarly, the goal of Fisher Criterion Score with SVM (Fisher-SVM) and Recursive Feature Elimination algorithm with SVM (RFE-SVM) is also restricting the number of features selected. They were introduced in \cite{GuyGunNikZad06,Guyon02}, respectively. Finally, $\nameM$ was also compared with the classical \mbox{(SVM-$\ell_1$)}.

To compare the classifiers, we used the recognized classification performance metrics: the accuracy (ACC) and the area under the curve (AUC).  The accuracy is computed as $$\text{ACC}=\dfrac{TP+TN}{TP+TN+FP+FN},$$ where TP are true positives, TN are true negatives, FP false positives and FN false negatives. The area under the curve is computed as $$\text{AUC}=\dfrac{\dfrac{TP}{TP+FN}+\dfrac{TN}{TN+FP}}{2}.$$

We followed the Ten Fold Cross Validation procedure (TFCV) to obtain these metrics. It consists in randomly partitioning the dataset into 10 subsets. In each iteration, nine of these subsets constitute the training set and the remaining partition is the test set. The performance of the classifier is evaluated by the independent test set.  

 Although many of the tested dataset are known to contain outliers, in order to check the robustness of the classifiers, we  perturbed the original data including label noise  and SVM outliers. The label noise can come from several real situations  which produce (intentionally or unintentionally) errors in the class of the individuals by adding outliers to the dataset, see \cite{Salgado2016}. For this purpose, we randomly changed the class of $5\%$  of the  individuals in the training sample. On the other hand, for adding SVM outliers, we solved problem \text{(SVM-$\ell_1$)} in the training sample. We sorted the individuals of each class in non-increasing order according to the result of $y_i\left(\sum_{k=1}^{d}\left(w_k^*x_{ik}\right)+b^*\right),$ where $i\in N, k\in D.$ Note that $w^*$ and $b^*$ are the optimal solution values of (SVM-$\ell_1$). We changed the class of the first $5\%$ ranked individuals.

The tables in which this information is depicted are structured as follows: in the first column the name of the classifier is shown and in the second column the value of parameter $B$ is reported. The following five columns describe information about the best case of $C$ when 5\% of the dataset contains label noise. The best case of $C$ is the value of parameter $C$ that reported the highest accuracy. Note that in the event of a tie, the best case is the one that provided a larger area under the curve. If this value is also the same, the best case is the one that took the least time. More concretely, the third column of the tables reports the above described value of parameter $C$, the next one shows the average time of the ten iterations, the following column depicts the average number of selected features for each classifier. Finally, the penultimate and last column of this group report the average ACC and AUC obtained for the ten iterations. The following five columns depict information about the best case of $C$ when 5\% of the dataset are SVM outliers. These columns are structured as the previous ones.

\begin{table}[p]
  \centering
  \resizebox{!}{0.23\textheight}{
    \begin{tabular}{lr|rrrrr|rrrrr}
    \hline
    \multicolumn{12}{l}{COLON, $n=62,d=2000$} \\
    \hline
          &       & \multicolumn{5}{c|}{  5\% label noise } & \multicolumn{5}{c}{  5\% SVM outliers } \\
    \hline
    Form. & \multicolumn{1}{l|}{$B$} & \multicolumn{1}{l}{$C$} & \multicolumn{1}{l}{Time} & \multicolumn{1}{l}{Av. F} & \multicolumn{1}{l}{ACC} & \multicolumn{1}{l|}{AUC} & \multicolumn{1}{l}{$C$} & \multicolumn{1}{l}{Time} & \multicolumn{1}{l}{Av. F} & \multicolumn{1}{l}{ACC} & \multicolumn{1}{l}{AUC} \\
    \hline
  \multirow{5}[2]{*}{\nameM} & 34    & 1     & 18.13 & 32.00    & 85.83\% & 83.75\% & 1     & 16.92 & 34.00 & 84.58\% & 83.75\% \\
          & 23    & 1     & 20.11 & 23.00    & 87.50\% & 85.00\% & 1     & 18.93 & 23.00 & 83.75\% & 82.50\% \\
          & 17    & 1     & 19.99 & 17.00   &87.50\% & 86.25\% & 1     & 18.86 & 17.00 & 87.08\% & 86.25\% \\
          & 11    & 1     & 18.42 & 11.00    & 84.58\% & 82.50\% & 1     & 14.08 & 11.00 & 87.08\% & 85.00\% \\
          & 7     & 1     & 16.53 & 7.00     & 82.50\% & 81.25\% & 1     & 14.70 & 7.00  & 83.75\% & 81.25\% \\
    \hline
    \multirow{5}[2]{*}{(FS-SVM)} & 34    & 1     & 0.80  & 32.60  & 85.83\% & 83.75\% & 10    & 3.33  & 34.00 & 87.08\% & 86.25\% \\
          & 23    & 1     & 3.00  & 23.00    & 87.08\% & 85.00\% & 10    & 184.12 & 23.00 & 85.42\% & 83.75\% \\
          & 17    & 1     & 33.63 & 17.00    & 87.50\% & 86.25\% & 100   & 524.55 & 17.00 & 83.75\% & 82.50\% \\
          & 11    & 10    & 249.24 & 11.00    & 80.83\% & 78.75\% & 1     & 180.45 & 11.00 & 85.42\% & 85.00\% \\
          & 7     & 1     & 20.11 & 7.00     & 85.83\% & 85.00\% & 10    & 885.16 & 7.00  & 83.75\% & 82.50\% \\
    \hline
        \multirow{5}[2]{*}{(Fisher-SVM)} & {34} & {0.1} & {0.01} & {34.00} & {82.08\%} & {81.25\%} & {0.1} & {0.01} & {34.00} & {85.42\%} & {86.25\%}\\
          & {23} & {1} &{0.01} &{23.00} & {82.08\%} & {81.25\%} & {0.1} &{0.01} &{23.00} & {85.42\%} & {86.25\%} \\
          & {17} & {1} & {0.01} & {17.00} & {85.42\%} & {85.00\%} & {1} & {0.01} & {17.00} & {86.25\%} & {86.25\%} \\
          & {11} & {10} &{0.01} & {11.00} &{85.83\%} & {86.25\%} & {0.1} &{0.01} &{11.00} & {83.75\%} & {81.25\%} \\
          & {7} &{\textbf{1}} {\textbf{0.01}} & {\textbf{7.00}} & {\textbf{88.75\%}} & {\textbf{88.75\%}} &{0.1} & {0.01} & {7.00} & {82.50\%} & {78.75\%} \\
    \hline
    \multirow{5}[2]{*}{{(RFE-SVM)}} & {34} & {0.1} & {0.01} & {34.00} & {85.42\%} & {84.17\%} & {1} & {0.01} & {34.00} & {83.75\%} & {83.75\%} \\
          & {23} & {0.1} & {0.01} & {23.00} & {85.42\%} & {82.92\%} & {1} & {0.01} & {23.00} & {82.08\%} & {82.50\%} \\
          & {17} & {10} & {0.01} & {17.00} & {83.75\%} & {82.92\%} & {1} & {0.01} & {17.00} & {82.08\%} & {82.50\%} \\
          & {11} & {0.1} & {0.01} & {11.00} & {82.08\%} & {78.75\%} & {1} & {0.01} & {11.00} & {82.08\%} & {82.50\%} \\
          & {7} & {1} & {0.01} & {7.00} & {82.08\%} & {80.83\%} & {1} & {0.01} & {7.00} & {85.42\%} & {85.00\%} \\
    \hline
    (RL-$\ell_1$-M) & -     & 1     & 23.29 & 35.80  & 87.08\% & 85.00\% & \textbf{1}     & \textbf{21.66} & \textbf{38.20} & \textbf{88.75\%} & \textbf{87.50\%} \\
    \hline
    (SVM-$\ell_1$) & -     & 1     & 0.67  & 32.70  & 85.83\% & 83.75\% & 1     & 0.68  & 36.70 & 84.17\% & 83.75\% \\
    \hline
    \end{tabular}%
    }
    \caption{Best average ACC and AUC for the Colon dataset.}
  \label{tab:TFCVColon}%
  
  \vspace{0.3cm}
   \resizebox{!}{0.23\textheight}{
 
    \begin{tabular}{lr|rrrrr|rrrrr}
    \hline
    \multicolumn{12}{l}{Leukemia, $n=72,d=5327$}\\
    \hline
          &       & \multicolumn{5}{c|}{  5\% label noise } & \multicolumn{5}{c}{  5\% SVM outliers }\\
    \hline
    Form. & \multicolumn{1}{l|}{$B$} & \multicolumn{1}{l}{$C$} & \multicolumn{1}{l}{Time} & \multicolumn{1}{l}{Av. F} & \multicolumn{1}{l}{ACC} & \multicolumn{1}{l|}{AUC} & \multicolumn{1}{l}{$C$} & \multicolumn{1}{l}{Time} & \multicolumn{1}{l}{Av. F} & \multicolumn{1}{l}{ACC} & \multicolumn{1}{l}{AUC} \\
    \hline
    \multirow{5}[2]{*}{\nameM} & 48    & 0.1   & 29.26 & 11.70  & 92.86\% & 93.75\% & 0.1   & 33.48 & 13.70 & 94.29\% & 95.00\% \\
          & 32    & 0.1   & 29.31 & 11.70  & 92.86\% & 93.75\% & \textbf{100}   & \textbf{120.56} & \textbf{32.0}0 & \textbf{94.60\%} & \textbf{95.54\%} \\
          & 24    & 0.1   & 29.75 & 11.70  & 92.86\% & 93.75\% & 0.1   & 33.45 & 12.60 & 94.29\% & 95.00\% \\
          & 16    & 0.1   & 29.50 & 11.80  & 92.86\% & 93.75\% & 0.1   & 33.48 & 12.60 & 94.29\% & 95.00\% \\
          & 10    & \textbf{0.1}   & \textbf{38.63} & \textbf{9.90 }  & \textbf{92.86\%} & \textbf{93.75\%} & 0.1   & 38.16 & 9.90  & 94.29\% & 95.00\% \\
    \hline
    \multirow{5}[2]{*}{(FS-SVM)} & 48    & 0.1   & 1.99  & 22.80  & 90.00\% & 90.83\% & 100   & 6.85  & 48.00 & 92.86\% & 93.33\% \\
          & 32    & 0.1   & 2.07  & 22.80  & 90.00\% & 90.83\% & 10    & 517.30 & 32.00 & 94.29\% & 94.58\% \\
          & 24    & 0.1   & 2.08  & 18.80  & 90.00\% & 90.83\% & 1     & 1402.60 & 24.00 & 93.17\% & 93.87\% \\
          & 16    & 100   & 956.59 & 16.00    & 90.32\% & 90.95\% & 1     & 1798.81 & 16.00 & 92.06\% & 93.15\% \\
          & 10    & 0.1   & 2.55  & 10.00    & 90.00\% & 90.83\% & 0.1   & 4.38  & 10.00 & 88.57\% & 89.17\%\\
    \hline
    \multirow{5}[2]{*}{{(Fisher-SVM)}} & {48} & {1} & {0.02} & {48.00} & {82.06\%} & {81.49\%} & {0.1} & {0.02} & {48.00} & {90.32\%} & {90.54\%} \\
          & {32} & {0.1} & {0.01} & {32.00} & {79.52\%} & {79.94\%} & {0.01} & {0.02} & {32.00} & {91.43\%} & {90.83\%} \\
          & {24} & {10} & {0.01} & {24.00} & {79.52\%} & {79.94\%} & {1} & {0.01} & {24.00} & {92.86\%} & {92.92\%} \\
          & {16} & {0.1} & {0.01} & {16.00} & {74.92\%} & {74.82\%} & {0.1} & {0.01} & {16.00} & {91.75\%} & {91.79\%} \\
          & {10} & {100} & {0.02} & {10.00} & {78.10\%} & {78.27\%} & {0.1} & {0.01} & {10.00} & {90.95\%} & {91.61\%} \\
    \hline
    \multirow{5}[2]{*}{{(RFE-SVM)}} & {48} & {0.01} & {0.01} & {48.00} & {87.46\%} & {88.04\%} & {0.01} & {0.01} & {48.00} & {90.00\%} & {90.00\%} \\
          & {32} & {0.01} & {0.01} & {32.00} & {86.03\%} & {86.79\%} & {0.01} & {0.01} & {32.00} & {92.86\%} & {92.92\%} \\
          & {24} & {0.01} & {0.01} & {24.00} & {87.46\%} & {88.04\%} & {0.01} & {0.01} & {24.00} & {94.29\%} & {94.58\%} \\
          & {16} & {0.01} & {0.01} & {16.00} & {86.03\%} & {86.79\%} & {0.01} & {0.01} & {16.00} & {92.86\%} & {92.92\%} \\
          & {10} & {0.01} & {0.01} & {10.00} & {87.46\%} & {88.04\%} & {0.01} & {0.01} & {10.00} & {92.86\%} & {93.33\%} \\
    \hline
    (RL-$\ell_1$-M) & -     & 0.1   & 258.50 & 11.80  & 91.75\% & 93.04\% & 0.1   & 632.51 & 12.00 & 94.29\% & 95.00\% \\
    \hline
    (SVM-$\ell_1$) & -     & 0.1   & 2.52  & 11.40  & 90.00\% & 90.83\% & 100   & 1.96  & 51.30 & 92.86\% & 93.33\% \\
    \hline
    \end{tabular}%
    }
     \caption{Best average ACC and AUC for the Leukemia dataset.}
  \label{tab:TFCVLeukemia}%

\end{table}%

\begin{table}[p]
  \centering
  \vspace{0.3cm}
 \resizebox{!}{0.23\textheight}{
    \begin{tabular}{lr|rrrrr|rrrrr}
    \hline
    \multicolumn{12}{l}{DLBCL, $n=77,d=7129$}\\
    \hline
          &       & \multicolumn{5}{c|}{  5\% label noise } & \multicolumn{5}{c}{  5\% SVM outliers } \\
    \hline
    Form. & \multicolumn{1}{l|}{$B$} & \multicolumn{1}{l}{$C$} & \multicolumn{1}{l}{Time} & \multicolumn{1}{l}{Av. F} & \multicolumn{1}{l}{ACC} & \multicolumn{1}{l|}{AUC} & \multicolumn{1}{l}{$C$} & \multicolumn{1}{l}{Time} & \multicolumn{1}{l}{Av. F} & \multicolumn{1}{l}{ACC} & \multicolumn{1}{l}{AUC} \\
    \hline
    \multirow{5}[2]{*}{\nameM} & 44    & 100   & 20.10 & 43.30  & 88.75\% & 82.50\% & 10    & 42.63 & 43.70 & 92.50\% & 90.00\% \\
          & 29    & 100   & 59.70 & 29.00    & 86.25\% & 79.17\% & 100   & 85.01 & 29.00 & 91.25\% & 89.17\% \\
          & 22    & 10    & 74.67 & 22.00    & 88.75\% & 84.17\% & 10    & 96.91 & 22.00 & 90.00\% & 86.67\% \\
          & 15    & 10    & 86.92 & 15.00    & 85.50\% & 76.67\% & \textbf{10}    & \textbf{126.05} & \textbf{15.00} & \textbf{95.00\%} & \textbf{95.00\%} \\
          & 9     & 100   & 159.07 & 9.00     & 85.50\% & 78.33\% & 100   & 425.14 & 9.00  & 92.50\% & 90.00\% \\
    \hline
    \multirow{5}[2]{*}{(FS-SVM)} & 44    & 1     & 4.72  & 43.60  & 88.75\% & 82.50\% & 1     & 7.20  & 43.80 & 92.50\% & 90.00\% \\
          & 29    & 1     & 127.71 & 29.00    & 87.50\% & 81.67\% & 1     & 339.38 & 29.00 & 91.25\% & 87.50\% \\
          & 22    & 100   & 599.93 & 22.00    & 88.75\% & 85.83\% & 1     & 1137.20 & 22.00 & 92.50\% & 90.00\% \\
          & 15    & 1     & 1348.37 & 15.00    & 85.50\% & 83.33\% & 10    & 1708.56 & 15.00 & 93.75\% & 92.50\% \\
          & 9     & 10    & 1494.80 & 9.00     & 86.75\% & 79.17\% & 1     & 1678.26 & 9.00  & 91.25\% & 85.83\% \\
           \hline
    \multirow{5}[2]{*}{{(Fisher-SVM)}} & {44} & {0.1} & {0.01} & {44.00} & {81.75\%} & {64.17\%} & {0.1} & {0.01} & {44.00} & {83.00\%} & {85.42\%} \\
          & {29} & {0.1} & {0.01} & {29.00} & {76.75\%} & {52.50\%} & {0.1} & {0.01} & {29.00} & {78.50\%} & {75.42\%} \\
          & {22} & {1} & {0.01} & {22.00} & {85.00\%} & {85.00\%} & {0.1} & {0.01} & {22.00} & {79.75\%} & {74.58\%} \\
          & {15} & {1} & {0.01} & {15.00} & {78.00\%} & {67.08\%} & {0.1} & {0.01} & {15.00} & {81.00\%} & {70.42\%} \\
          & {9} & {1} & {0.02} & {9.00} & {80.00\%} & {66.67\%} & {0.1} & {0.01} & {9.00} & {78.00\%} & {58.33\%} \\
    \hline
    \multirow{5}[2]{*}{{(RFE-SVM)}} & {44} & {1} & {0.02} & {44.00} & {83.75\%} & {81.67\%} & {0.1} & {0.01} & {44.00} & {91.25\%} & {90.83\%} \\
          & {29} & {1} & {0.02} & {29.00} & {81.25\%} & {81.67\%} & {0.1} & {0.01} & {29.00} & {91.25\%} & {90.83\%} \\
          & {22} & {0.1} & {0.01} & {22.00} & {82.50\%} & {80.67\%} & {0.1} & {0.01} & {22.00} & {91.25\%} & {89.17\%} \\
          & {15} & {0.1} & {0.01} & {15.00} & {82.50\%} & {77.33\%} & {0.1} & {0.01} & {15.00} & {95.00\%} & {93.33\%} \\
          & {9} & {10} & {0.01} & {9.00} & {80.00\%} & {80.83\%} & {0.1} & {0.01} & {9.00} & {91.25\%} & {87.50\%} \\
    \hline    (RL-$\ell_1$-M) & -     & \textbf{1}     & \textbf{74.40} & \textbf{51.10}  & \textbf{91.25\%} & \textbf{85.83\%} & 1     & 74.99 & 51.60 & 92.50\% & 91.67\% \\
    \hline
    (SVM-$\ell_1$) & -     & 1     & 3.86  & 43.80  & 88.75\% & 82.50\% & 1     & 3.59  & 47.20 & 92.50\% & 90.00\% \\
    \hline
    \end{tabular}%
    }
    \caption{Best average ACC and AUC for the DLBCL dataset.}
  \label{tab:TFCVDLBCL}%
    \vspace{0.3cm}
   \resizebox{!}{0.23\textheight}{

    \begin{tabular}{lr|rrrrr|rrrrr}
    \hline
      \multicolumn{12}{l}{SONAR, $n=208,d=60$} \\
    \hline
          &       & \multicolumn{5}{c|}{  5\% label noise } & \multicolumn{5}{c}{  5\% SVM outliers } \\
    \hline
    Form. & \multicolumn{1}{l|}{$B$} & \multicolumn{1}{l}{$C$} & \multicolumn{1}{l}{Time} & \multicolumn{1}{l}{Av. F} & \multicolumn{1}{l}{ACC} & \multicolumn{1}{l|}{AUC} & \multicolumn{1}{l}{$C$} & \multicolumn{1}{l}{Time} & \multicolumn{1}{l}{Av. F} & \multicolumn{1}{l}{ACC} & \multicolumn{1}{l}{AUC} \\
    \hline
    \multirow{5}[2]{*}{\nameM} & 57    & 1     & 47.40 & 32.70  & 76.42\% & 76.05\% & 1     & 81.06 & 33.10 & 75.94\% & 75.59\% \\
          & 38    & \textbf{10}    & \textbf{12.16} & \textbf{38.00}    & \textbf{77.52\%} & \textbf{77.03\%} & 0.1   & 237.67 & 7.50  & 75.51\% & 75.01\% \\
          & 29    & 1     & 63.37 & 28.90  & 76.47\% & 75.97\% & 10    & 33.59 & 29.00 & 75.94\% & 75.48\% \\
          & 19    & 10    & 148.01 & 19.00    & 75.94\% & 75.09\% & \textbf{10}    & \textbf{220.12} & \textbf{19.00} & \textbf{76.94\%} & \textbf{76.90\%} \\
          & 11    & 0.1   & 268.60 & 11.00    & 76.89\% & 76.76\% & 0.1   & 402.95 & 8.10  & 75.51\% & 75.01\% \\
    \hline
    \multirow{5}[2]{*}{(FS-SVM)} & 57    & 0.1   & 0.07  & 14.40  & 76.99\% & 76.73\% & 0.1   & 0.08  & 17.80 & 72.61\% & 72.19\% \\
          & 38    & 0.1   & 0.06  & 14.40  & 76.99\% & 76.73\% & 0.1   & 0.07  & 17.80  & 72.61\% & 72.19\% \\
          & 29    & 0.1   & 0.06  & 14.40  & 76.99\% & 76.73\% & 10    & 678.46 & 29.00 & 73.61\% & 73.81\% \\
          & 19    & 0.1   & 0.07  & 14.40  & 76.99\% & 76.73\% & 10    & 753.49 & 19.00 & 74.09\% & 74.08\% \\
          & 11    & 0.1   & 0.07  & 8.20   & 76.99\% & 76.73\% & 1     & 137.99 & 11.00 & 73.18\% & 72.95\% \\
    \hline
    \multirow{5}[2]{*}{{(Fisher-SVM)}} & {57} & {1} & {0.03} & {57.00} & {75.94\%} & {76.76\%} & {0.1} & {0.03} & {57.00} & {73.61\%} & {74.17\%}]\\
          & {38} & {1} & {0.02} & {38.00} & {74.51\%} & {75.21\%} & {0.1} & {0.02} & {38.00} & {72.61\%} & {73.26\%} \\
          & {29} & {1} & {0.02} & {29.00} & {76.52\%} & {77.00\%} & {1} & {0.02} & {29.00} & {72.03\%} & {72.42\%} \\
          & {19} & {1} & {0.02} & {19.00} & {71.60\%} & {72.00\%} & {1} & {0.02} & {19.00} & {70.45\%} & {71.12\%} \\
          & {11} & {1} & {0.03} & {11.00} & {73.03\%} & {73.45\%} & {1} & {0.03} & {11.00} & {70.03\%} & {70.54\%} \\
    \hline
    \multirow{5}[2]{*}{{(RFE-SVM)}} & {57} & {1} & {0.03} & {57.00} & {75.94\%} & {76.76\%} & {0.1} & {0.03} & {57.00} & {73.61\%} & {74.26\%} \\
          & {38} & {1} & {0.02} & {38.00} & {71.18\%} & {71.76\%} & {0.01} & {0.02} & {38.00} & {74.09\%} & {74.49\%} \\
          & {29} & {0.1} & {0.02} & {29.00} & {72.13\%} & {72.66\%} & {10} & {0.02} & {29.00} & {75.04\%} & {75.44\%} \\
          & {19} & {1} & {0.02} & {19.00} & {73.61\%} & {73.94\%} & {0.1} & {0.02} & {19.00} & {71.65\%} & {72.26\%} \\
          & {11} & {0.1} & {0.03} & {11.00} & {72.13\%} & {72.44\%} & {0.1} & {0.03} & {11.00} & {71.18\%} & {71.66\%} \\
          \hline
    (RL-$\ell_1$-M) & -     & 1     & 1806.60 & 32.90  & 76.94\% & 76.47\% &      1     & 1805.03 & 33.10  & 76.37\% & 75.79\%\\
    \hline
    (SVM-$\ell_1$) & -     & 0.1   & 0.02  & 7.20   & 76.99\% & 76.73\% & 0.1   & 0.02  & 8.90  & 72.61\% & 72.19\% \\
    \hline
    \end{tabular}%
    }
     \caption{Best average ACC and AUC for the SONAR dataset.}
  \label{tab:TFCVSonar}%
 \end{table}%
In Table \ref{tab:TFCVColon}, the computational experiments on the Colon dataset are depicted. We can observe that the best classification performance when 5\% of the dataset contains label noise
is obtained by (Fisher-SVM), with the following metrics: 88.75\% of ACC and 88.75\% of AUC.  
On the other hand, the best classifier is \text{(RL-$\ell_1$-M)} when 5\% of the dataset are SVM outliers, providing 88.75\% ACC and 87.50\% AUC. Slightly smaller classification metrics are obtained by $\nameM$ selecting only eleven features on average in contrast with the 38.20 selected by \text{(RL-$\ell_1$-M)}.

The classification performance for the different formulations in the Leukemia dataset is depicted in Table \ref{tab:TFCVLeukemia}. Note that the best classifier when 5\% of the dataset contains label noise 
is $\nameM$, providing 92.86\% ACC and 93.75\% AUC. Moreover, observe that our model is able to classify a new individual by analyzing only 9.9 features on average. Similarly, the largest classification metrics (94.60\% ACC and 95.54\% AUC) are obtained by $\nameM$ when 5\% of the dataset are SVM outliers. 

The results from the DLBCL dataset are reported in Table \ref{tab:TFCVDLBCL}. As can be appreciated, when 5\% of the dataset contains label noise 
the best classifier is \mbox{(RL-$\ell_1$-M)}, providing 91.25\% ACC and 85.83\% AUC. However, this classifier used 51.1 features on average. Slightly smaller classification metrics are obtained by $\nameM$ and \text{(FS-SVM)} selecting fewer than half of the features on average. On the other hand, when 5\% of the dataset are SVM outliers, $\nameM$ is the classifier that provides the largest ACC (95\%) and AUC (95\%) with nearly a point and a half difference in percentage compared to the rest of the models. Observe that this classifier selected only fifteen features on average from the 7129 provided by the data.

\begin{table}[p]
  \centering

\vspace{0.3cm}
 \resizebox{!}{0.23\textheight}{
    \begin{tabular}{lr|rrrrr|rrrrr}
    \hline
        \multicolumn{12}{l}{IONO, $n=351,d=33$}\\
    \hline
          &       & \multicolumn{5}{c|}{  5\% label noise } & \multicolumn{5}{c}{  5\% SVM outliers } \\
    \hline
    Form. & \multicolumn{1}{l|}{$B$} & \multicolumn{1}{l}{$C$} & \multicolumn{1}{l}{Time} & \multicolumn{1}{l}{Av. F} & \multicolumn{1}{l}{ACC} & \multicolumn{1}{l|}{AUC} & \multicolumn{1}{l}{$C$} & \multicolumn{1}{l}{Time} & \multicolumn{1}{l}{Av. F} & \multicolumn{1}{l}{ACC} & \multicolumn{1}{l}{AUC} \\
    \hline
     \multirow{5}[2]{*}{\nameM}  & 33    & 1     & 75.71 & 25.00    & 88.06\% & 84.34\% & 1     & 108.15 & 25.30 & 86.57\% & 82.58\% \\
          & 22    & 1     & 119.80 & 21.90  & 87.96\% & 84.00\% & 1     & 111.09 & 22.00 & 85.19\% & 80.49\% \\
          & 17    & 10    & 34.13 & 17.00    & 87.13\% & 83.35\% & 1     & 134.82 & 17.00 & 86.02\% & 81.98\% \\
          & 11    & 10    & 185.26 & 11.00    & 86.02\% & 81.64\% & 1     & 206.58 & 11.00 & 86.57\% & 82.08\% \\
          & 7     & 10    & 197.31 & 7.00     & 86.76\% & 83.74\% & 1     & 367.34 & 7.00  & 85.37\% & 80.87\% \\
    \hline
    \multirow{5}[2]{*}{(FS-SVM)} & 33    & 1     & 0.07  & 33.00    & 87.04\% & 82.95\% & 100   & 0.06  & 33.00 & 78.61\% & 71.63\% \\
          & 22    & 100   & 0.92  & 22.00    & 87.22\% & 84.13\% & 100   & 3.05  & 22.00 & 72.59\% & 62.55\% \\
          & 17    & 1     & 1.89  & 17.00    & 87.41\% & 83.90\% & 100   & 16.82 & 17.00 & 72.04\% & 62.28\% \\
          & 11    & \textbf{1}     & \textbf{3.81}  & \textbf{11.00 }   & \textbf{88.15\%} & \textbf{84.66\%} & 100   & 40.39 & 11.00 & 73.70\% & 63.58\% \\
          & 7     & 100   & 20.20 & 7.00     & 87.69\% & 84.23\% & 100   & 35.78 & 7.00  & 71.94\% & 61.06\% \\
          \hline
          \multirow{5}[2]{*}{{(Fisher-SVM)}} & {33} & {0.1} & {0.04} & {33.00} & {87.96\%} & {84.00\%} & {0.01} & {0.04} & {33.00} & {86.85\%} & {82.46\%} \\
          & {22} & {0.01} & {0.04} & {22.00} & {85.19\%} & {80.65\%} & {0.1} & {0.04} & {22.00} & {86.94\%} & {82.13\%} \\
          & {17} & {0.01} & {0.04} & {17.00} & {84.35\%} & {79.83\%} & {1} & {0.04} & {17.00} & {86.11\%} & {81.64\%} \\
          & {11} & {1} & {0.06} & {11.00} & {84.35\%} & {80.17\%} & {\textbf{100}} & {\textbf{0.05}} & {\textbf{11.00}} & {\textbf{87.78\%}} & {\textbf{83.45\%}} \\
          & {7} & {0.01} & {0.05} & {7.00} & {81.85\%} & {76.37\%} & {10} & {0.05} & {7.00} & {87.50\%} & {82.90\%} \\
    \hline
    \multirow{5}[2]{*}{{(RFE-SVM)}} & {33} & {0.1} & {0.05} & {33.00} & {84.17\%} & {80.36\%} & {0.01} & {0.04} & {33.00} & {86.85\%} & {82.46\%} \\
          & {22} & {0.1} & {0.04} & {22.00} & {84.44\%} & {80.61\%} & {0.01} & {0.03} & {22.00} & {86.85\%} & {82.46\%} \\
          & {17} & {1} & {0.05} & {17.00} & {84.44\%} & {80.24\%} & {0.01} & {0.04} & {17.00} & {86.94\%} & {82.63\%} \\
          & {11} & {1} & {0.05} & {11.00} & {83.06\%} & {78.57\%} & {0.01} & {0.05} & {11.00} & {86.57\%} & {81.91\%} \\
          & {7} & {1} & {0.06} & {7.00} & {82.87\%} & {78.42\%} & {0.01} & {0.04} & {7.00} & {85.37\%} & {80.53\%} \\
    \hline
    (RL-$\ell_1$-M) & -     & 1     & 1805.00 & 25.60  & 88.06\% & 84.34\% &    1   &   1804.40    &   24.60    &   85.74\% &  81.59\%\\
    \hline
    (SVM-$\ell_1$) & -     & 1     & 0.01  & 28.00    & 87.04\% & 82.95\% & 1     & 0.01  & 26.40 & 86.02\% & 81.81\%\\
    \hline
    \end{tabular}%
    }
    \caption{Best average ACC and AUC for the IONO dataset.}
  \label{tab:TFCVIONO}%
    \vspace{0.3cm}
     \resizebox{!}{0.23\textheight}{
     \begin{tabular}{lr|rrrrr|rrrrr}
    \hline
    \multicolumn{12}{l}{Arrhythmia, $n=420,d=258$} \\
    \hline
          &       & \multicolumn{5}{c|}{  5\% label noise } & \multicolumn{5}{c}{  5\% SVM outliers }\\
    \hline
    Form. & \multicolumn{1}{l|}{$B$} & \multicolumn{1}{l}{$C$} & \multicolumn{1}{l}{Time} & \multicolumn{1}{l}{Av. F} & \multicolumn{1}{l}{ACC} & \multicolumn{1}{l|}{AUC} & \multicolumn{1}{l}{$C$} & \multicolumn{1}{l}{Time} & \multicolumn{1}{l}{Av. F} & \multicolumn{1}{l}{ACC} & \multicolumn{1}{l}{AUC} \\
    \hline
  \multirow{5}[2]{*}{\nameM} & 81    & 1     & 1629.88 & 79.60  & 73.57\% & 72.08\% & 1     & 1866.23 & 77.30 & 73.33\% & 71.62\% \\
          & 54    & 1     & 1793.33 & 54.00    & 76.43\% & 75.17\% & 1     & 1860.52 & 54.00 & 70.95\% & 69.24\% \\
          & 32    & 10    & 1834.59 & 32.00    & 75.48\% & 74.68\% & 1     & 1830.04 & 32.00 & 72.38\% & 71.48\% \\
          & 16    & 1     & 1832.23 & 16.00    & 75.48\% & 74.01\% & 1     & 1827.54 & 16.00 & 74.52\% & 72.63\% \\
          & 8     & 1     & 1845.67 & 8.00     & 72.14\% & 70.31\% & 10    & 1846.35 & 8.00  & 70.71\% & 68.76\% \\
    \hline
    \multicolumn{1}{c}{\multirow{5}[2]{*}{(FS-SVM)}} & 81    & 1     & 267.70 & 81.00    & 75.00\% & 73.71\% & 1     & 615.29 & 80.80 & 70.95\% & 69.68\% \\
          & 54    & 10    & 919.43 & 54.00    & 73.57\% & 73.18\% & 1     & 890.52 & 54.00 & 72.62\% & 70.49\% \\
          & 32    & 10    & 882.25 & 32.00    & 72.86\% & 72.00\% & 10    & 755.65 & 32.00 & 73.33\% & 72.03\% \\
          & 16    & 10    & 862.04 & 16.00    & 76.19\% & 74.55\% & 10    & 805.40 & 16.00 & 75.95\% & 74.13\% \\
          & 8     & \textbf{10}   & \textbf{745.59} & \textbf{8.00}     & \textbf{76.67\%} & \textbf{74.74\%} & 100   & 875.68 & 8.00  & 73.57\% & 71.46\%\\
    \hline
    \multirow{5}[2]{*}{{(Fisher-SVM)}} & {81} & {0.1} & {0.05} & {81.00} & {72.38\%} & {69.00\%} & {10} & {0.06} & {81.00} & {71.67\%} & {67.99\%}\\
          & {54} & {1} & {0.05} & {54.00} & {70.95\%} & {67.39\%} & {10} & {0.06} & {54.00} & {70.95\%} & {67.65\%} \\
          & {32} & {10} & {0.02} & {32.00} & {68.81\%} & {64.65\%} & {10} & {0.04} & {32.00} & {69.29\%} & {65.35\%} \\
          & {16} & {100} & {0.02} & {16.00} & {66.67\%} & {62.52\%} & {10} & {0.04} & {16.00} & {67.62\%} & {63.76\%} \\
          & {8} & {100} & {0.02} & {8.00} & {62.38\%} & {57.32\%} & {10} & {0.08} & {8.00} & {66.90\%} & {63.76\%} \\
    \hline
    \multirow{5}[2]{*}{{(RFE-SVM)}} & {81} & {0.01} & {0.06} & {81.00} & {70.24\%} & {67.36\%} & {1} & {0.07} & {81.00} & {70.71\%} & {67.81\%} \\
          & {54} & {0.1} & {0.06} & {54.00} & {69.76\%} & {66.46\%} & {0.01} & {0.05} & {54.00} & {72.38\%} & {69.34\%} \\
          & {32} & {0.1} & {0.04} & {32.00} & {70.71\%} & {67.64\%} & {1} & {0.04} & {32.00} & {73.33\%} & {70.63\%} \\
          & {16} & {0.1} & {0.03} & {16.00} & {70.48\%} & {67.43\%} & {0.1} & {0.02} & {16.00} & {69.76\%} & {65.86\%} \\
          & {8} & {0.1} & {0.04} & {8.00} & {68.33\%} & {65.00\%} & {0.1} & {0.03} & {8.00} & {69.52\%} & {65.36\%} \\
    \hline
    (RL-$\ell_1$-M) & -     & 1     & 1802.37 & 84.70 & 73.57\% & 72.07\% &      \textbf{1}     & \textbf{1801.78} & \textbf{79.30}  & \textbf{75.95\%} & \textbf{74.56\%}  \\
    \hline
    (SVM-$\ell_1$) & -     & 1     & 0.25  & 97.00    & 74.52\% & 73.15\% & 1     & 0.27  & 103.40 & 71.90\% & 69.93\% \\
    \hline
    \end{tabular}%
    }
    \caption{Best average ACC and AUC for the Arrhythmia dataset.}
  \label{tab:TFCVArrhythmia}%
  
\end{table}%

The computational results for the Sonar dataset are shown in Table \ref{tab:TFCVSonar}. When 5\% of the dataset contains label noise, 
 the largest ACC (77.52\%) and AUC (77.03\%) are obtained by $\nameM$. Likewise, in the case of having a 5\% of SVM outliers, the best classification performance (76.94\% ACC and 76.90\% AUC) is reported by $\nameM$. 

In Table \ref{tab:TFCVIONO}, the classification performance of the analyzed formulations on the Iono dataset is provided. Although the highest performance metrics when 5\% of the dataset contains label noise 
are reported by $\text{(FS-SVM)}$ whose average ACC is $88.15\%$ and average AUC is $84.66\%$, slightly smaller accuracy and area under the curve (88.06\% and 84.34\% respectively) are provided by $\nameM$ and $\text{(RL-$\ell_1$-M)}$. On the contrary, when 5\% of the dataset are SVM outliers, the best classifier  is (Fisher-SVM).

The experiments carried out on the Arrhythmia dataset are depicted in Table \ref{tab:TFCVArrhythmia}. Observe that when 5\% of the dataset contains label noise,
$\text{(FS-SVM)}$ is the classifier that provides the largest ACC (76.67\%) and AUC (74.74\%) on average. However, when 5\% of the dataset are SVM outliers, $\text{(RL-$\ell_1$-M)}$ reported the largest classification metrics: 75.95\% ACC and 74.56\% AUC. Slightly smaller classification metrics are obtained by $\text{(FS-SVM)}$ selecting only 16 features on average instead of 79.30.

\begin{table}[p]
  \centering
 
   \resizebox{!}{0.23\textheight}{

    \begin{tabular}{lr|rrrrr|rrrrr}
    \hline
         \multicolumn{12}{l}{Wdbc, $n=569,d=30$} \\
    \hline
          &       & \multicolumn{5}{c|}{  5\% label noise } & \multicolumn{5}{c}{  5\% SVM outliers } \\
    \hline
    Form. & \multicolumn{1}{l|}{$B$} & \multicolumn{1}{l}{$C$} & \multicolumn{1}{l}{Time} & \multicolumn{1}{l}{Av. F} & \multicolumn{1}{l}{ACC} & \multicolumn{1}{l|}{AUC} & \multicolumn{1}{l}{$C$} & \multicolumn{1}{l}{Time} & \multicolumn{1}{l}{Av. F} & \multicolumn{1}{l}{ACC} & \multicolumn{1}{l}{AUC} \\
    \hline
    \multirow{5}[2]{*}{\nameM}& 29    & \textbf{10}    & \textbf{34.93} & \textbf{16.40}  & \textbf{97.54\%} & \textbf{97.06\%} & 1     & 39.05 & 9.90  & 97.54\% & 97.06\% \\
          & 19    & 10    & 35.41 & 16.40  & 97.37\% & 96.83\% & \textbf{1}     & \textbf{39.00} & \textbf{9.90}  & \textbf{97.54\%} & \textbf{97.06\%} \\
          & 15    & 10    & 35.15 & 14.90  & 97.19\% & 96.69\% & 1     & 39.17 & 8.80  & 97.54\% & 97.06\% \\
          & 10    & 1     & 50.45 & 9.80   & 96.67\% & 95.97\% & 1     & 39.93 & 8.90  & 97.54\% & 97.06\% \\
          & 6     & 1     & 62.25 & 6.00     & 96.49\% & 95.63\% & 1     & 44.39 & 6.00  & 97.19\% & 96.69\% \\
    \hline
    \multirow{5}[2]{*}{(FS-SVM)} & 29    & 1     & 0.10  & 22.40  & 96.49\% & 95.56\% & 1     & 0.12  & 22.30 & 96.13\% & 95.00\% \\
          & 19    & 1     & 0.10  & 13.80  & 96.49\% & 95.56\% & 1     & 0.12  & 14.30 & 96.13\% & 95.00\% \\
          & 15    & 1     & 0.11  & 12.20  & 96.49\% & 95.56\% & 1     & 0.14  & 13.90 & 96.13\% & 95.00\% \\
          & 10    & 1     & 0.18  & 10.00    & 96.49\% & 95.56\% & 1     & 0.36  & 10.00 & 95.96\% & 94.78\% \\
          & 6     & 1     & 0.25  & 6.00     & 96.66\% & 95.79\% & 10    & 5.45  & 6.00  & 96.67\% & 95.97\% \\
    \hline
    \multirow{5}[2]{*}{{(Fisher-SVM)}} & {29} & {0.1} & {0.13} & {29.00} & {94.91\%} & {95.35\%} & {1} & {0.27} & {29.00} & {95.60\%} & {95.70\%} \\
          & {19} & {1} & {0.19} & {19.00} & {93.14\%} & {93.65\%} & {0.1} & {0.28} & {19.00} & {95.43\%} & {95.70\%} \\
          & {15} & {10} & {0.20} & {15.00} & {93.50\%} & {93.85\%} & {0.01} & {0.17} & {15.00} & {93.85\%} & {92.62\%} \\
          & {10} & {10} & {0.20} & {10.00} & {94.19\%} & {94.59\%} & {1} & {0.26} & {10.00} & {93.50\%} & {93.40\%} \\
          & {6} & {0.1} & {0.17} & {6.00} & {92.09\%} & {92.18\%} & {0.1} & {0.18} & {6.00} & {92.44\%} & {90.75\%} \\
    \hline
    \multirow{5}[2]{*}{{(RFE-SVM)}} & {29} & {10} & {0.14} & {29.00} & {90.86\%} & {90.71\%} & {1} & {0.45} & {29.00} & {95.78\%} & {95.85\%} \\
          & {19} & {10} & {0.18} & {19.00} & {91.22\%} & {90.92\%} & {1} & {0.39} & {19.00} & {95.78\%} & {95.95\%} \\
          & {15} & {10} & {0.18} & {15.00} & {91.57\%} & {91.20\%} & {1} & {0.34} & {15.00} & {95.43\%} & {95.57\%} \\
          & {10} & {10} & {0.19} & {10.00} & {90.69\%} & {90.15\%} & {1} & {0.23} & {10.00} & {95.61\%} & {95.82\%} \\
          & {6} & {0.1} & {0.14} & {6.00} & {90.34\%} & {90.03\%} & {0.1} & {0.22} & {6.00} & {95.61\%} & {95.71\%} \\
    \hline
    (RL-$\ell_1$-M) & -     & 1     & 1807.89 & 9.40   & 96.67\% & 95.97\% &     1     & 1111.60 & 8.70   & 97.54\% & 97.06\% \\
    \hline
    (SVM-$\ell_1$) & -     & 1     & 0.02  & 11.20  & 96.49\% & 95.56\% & 1     & 0.02  & 13.30 & 96.13\% & 95.00\% \\
    \hline
    \end{tabular}%
    }
     \caption{Best average ACC and AUC for the Wdbc dataset.}
  \label{tab:TFCVWdbc}%

\vspace{0.3cm}
 \resizebox{!}{0.23\textheight}{
    \begin{tabular}{lr|rrrrr|rrrrr}
    \hline
       \multicolumn{12}{l}{Mfeat, $n=2000,d=649$} \\
    \hline
          &       & \multicolumn{5}{c|}{  5\% label noise } & \multicolumn{5}{c}{  5\% SVM outliers } \\
    \hline
    Form. & \multicolumn{1}{l|}{$B$} & \multicolumn{1}{l}{$C$} & \multicolumn{1}{l}{Time} & \multicolumn{1}{l}{Av. F} & \multicolumn{1}{l}{ACC} & \multicolumn{1}{l|}{AUC} & \multicolumn{1}{l}{$C$} & \multicolumn{1}{l}{Time} & \multicolumn{1}{l}{Av. F} & \multicolumn{1}{l}{ACC} & \multicolumn{1}{l}{AUC} \\
    \hline
    \multirow{5}[2]{*}{\nameM} 
  & 202   & 0.1   & 516.33 & 35.50  & 99.90\% & 99.50\% &     0.1   & 1080.72 & 57.80  & 99.85\% & 99.25\%  \\
          & 135   & 0.1   & 516.10 & 35.50  & 99.90\% & 99.50\% &  0.1   & 1083.40 & 45.60  & 99.85\% & 99.25\% \\
          & 81    & 0.1   & 483.46 & 35.20  & 99.90\% & 99.50\% &     0.1   & 1092.25 & 45.00    & 99.90\% & 99.50\%  \\
          & 30    & 0.1   & 381.24 & 29.80  & 99.90\% & 99.50\% &   {\bf  0.1}   & {\bf 581.82} & {\bf 29.90}  & {\bf 99.90\%} & {\bf 99.50\%}   \\
          & 10    & 0.1   & 268.69 & 10.00   & 99.95\% & 99.75\% &       1     & 350.67 & 5.10   & 99.85\% & 99.25\% \\
    \hline
    \multicolumn{1}{c}{\multirow{5}[2]{*}{(FS-SVM)}} & 202   & 0.1   & 9.69  & 171.50 & 99.90\% & 99.50\% & 0.1   & 9.79  & 183.40 & 99.80\% & 99.22\%\\
          & 135   & 0.1   & 9.74  & 91.70  & 99.90\% & 99.50\% & 0.1   & 9.79  & 92.70 & 99.80\% & 99.22\% \\
          & 81    & 0.1   & 154.87 & 81.00    & 99.90\% & 99.50\% & 0.1   & 354.21 & 81.00 & 99.80\% & 99.22\% \\
          & 30    & \textbf{1}     & \textbf{997.05} & \textbf{30.00}    & \textbf{100.00\%} & \textbf{100.00\%} & 0.1   & 832.68 & 30.00 & 99.75\% & 98.97\% \\
          & 10    & 0.1   & 750.47 & 10.00    & 99.90\% & 99.72\% & 0.1   & 807.60 & 10.00 & 99.65\% & 98.69\%\\
    \hline
    \multirow{5}[2]{*}{{(Fisher-SVM)}} & {202} & {0.1} & {1.58} & {202.00} & {100.00\%} & {100.00\%} & {0.01} & {1.57} & {202.00} & {99.80\%} & {99.67\%} \\
          & {135} & {0.1} & {1.37} & {135.00} & {99.00\%} & {99.22\%} & {0.01} & {1.78} & {135.00} & {99.80\%} & {99.50\%} \\
          & {81} & {0.01} & {0.90} & {81.00} & {95.75\%} & {97.64\%} & {0.01} & {1.97} & {81.00} & {99.85\%} & {99.69\%} \\
          & {30} & {0.01} & {0.61} & {30.00} & {94.05\%} & {70.47\%} & {0.1} & {0.72} & {30.00} & {98.60\%} & {99.22\%} \\
          & {10} & {0.1} & {2.59} & {10.00} & {94.70\%} & {88.17\%} & {0.1} & {2.91} & {10.00} & {95.65\%} & {97.36\%} \\
    \hline
    \multirow{5}[2]{*}{{(RFE-SVM)}} & {202} & {0.1} & {1.08} & {202.00} & {94.90\%} & {83.43\%} & {0.01} & {2.32} & {202.00} & {97.30\%} & {98.28\%} \\
          & {135} & {0.01} & {0.75} & {135.00} & {94.90\%} & {83.58\%} & {0.01} & {1.23} & {135.00} & {98.05\%} & {98.69\%} \\
          & {81} & {0.01} & {0.47} & {81.00} & {94.85\%} & {83.40\%} & {0.01} & {0.93} & {81.00} & {98.15\%} & {98.75\%} \\
          & {30} & {0.01} & {0.30} & {30.00} & {94.90\%} & {83.58\%} & {0.01} & {0.32} & {30.00} & {99.40\%} & {99.44\%} \\
          & {10} & {0.01} & {0.81} & {10.00} & {94.20\%} & {83.32\%} & {0.01} & {1.49} & {10.00} & {98.60\%} & {99.00\%} \\
    \hline
    (RL-$\ell_1$-M) & -     & 0.1   & 1809.61 & 47.00    & 99.90\% & 99.50\% &   0.1
    &    1809.92   &  70.30     &    99.70\%   & 99.39\% \\
    \hline
    (SVM-$\ell_1$) & -     & 0.1   & 15.49 & 90.80  & 99.90\% & 99.50\% & 0.1   & 13.41 & 91.70 & 99.80\% & 99.22\% \\
    \hline
    \end{tabular}%
    }
    \caption{Best average ACC and AUC for the Mfeat dataset.}
  \label{tab:TFCVMfeat}%
\end{table}%

\begin{table}[ht]
    \resizebox{\textwidth}{!}{
     \begin{tabular}{lr|rrrrr|rrrrr}
    \hline
        \multicolumn{12}{l}{Lepiota, $n=8124,d=109$}\\
    \hline
          &       & \multicolumn{5}{c|}{  5\% label noise } & \multicolumn{5}{c}{  5\% SVM outliers } \\
    \hline
    Form. & \multicolumn{1}{l|}{$B$} & \multicolumn{1}{l}{$C$} & \multicolumn{1}{l}{Time} & \multicolumn{1}{l}{Av. F} & \multicolumn{1}{l}{ACC} & \multicolumn{1}{l|}{AUC} & \multicolumn{1}{l}{$C$} & \multicolumn{1}{l}{Time} & \multicolumn{1}{l}{Av. F} & \multicolumn{1}{l}{ACC} & \multicolumn{1}{l}{AUC} \\
    \hline
    \multirow{5}[2]{*}{\nameM}& 33    & 1     & 496.16 & 12.20  & 100.00\% & 100.00\% &      10    & 1741.94 &23.60  & 96.75\% & 96.64\%   \\
          & 22    & \textbf{1}     & \textbf{529.65} & \textbf{10.70}  & \textbf{100.00\%} & \textbf{100.00\%} &   10    & 1759.38 & 21.80  & 96.71\% & 96.60\%   \\
          & 17    & 1     & 508.80 & 11.70  & 99.99\% & 99.99\% &     10    & 1730.30 & 16.70 & 96.86\% & 96.75\%   \\
          & 11    & 1     & 841.15 & 11.00    & 99.95\% & 99.95\% &      0.1   & 370.89 & 8.60   & 97.51\% & 97.42\% \\
          & 7     & 10    & 962.17 & 7.00     & 99.84\% & 99.83\% &      \textbf{0.1}   & \textbf{414.28} & \textbf{6.40}   & \textbf{99.73\%} & \textbf{99.72\%}\\
    \hline
    \multirow{5}[2]{*}{(FS-SVM)} & 33    & 10    & 11.21 & 19.20  & 99.95\% & 99.95\% & 0.01  & 9.23  & 7.50  & 88.50\% & 63.81\% \\
          & 22    & 10    & 11.10 & 17.60  & 99.95\% & 99.95\% & 0.01  & 9.23  & 7.50  & 88.50\% & 63.81\% \\
          & 17    & 10    & 11.81 & 16.30  & 99.95\% & 99.95\% & 0.01  & 9.22  & 7.50  & 88.50\% & 63.81\% \\
          & 11    & 10    & 12.84 & 11.00    & 99.91\% & 99.91\% & 100   & 10.77 & 8.90  & 90.98\% & 89.22\% \\
          & 7     & 0.1   & 10.69 & 6.00     & 99.70\% & 99.69\% & 100   & 10.87 & 6.40  & 87.27\% & 86.77\% \\
    \hline
    \multirow{5}[2]{*}{{(Fisher-SVM)}} & {33} & {100} & {4.41} & {33.00} & {80.92\%} & {81.44\%} & {0.1} & {4.66} & {33.00} & {95.88\%} & {95.72\%} \\
          & {22} & {100} & {4.43} & {22.00} & {61.10\%} & {59.89\%} & {0.1} & {4.65} & {22.00} & {95.26\%} & {95.08\%} \\
          & {17} & {1} & {4.45} & {17.00} & {57.80\%} & {56.55\%} & {0.01} & {4.59} & {17.00} & {95.37\%} & {95.20\%} \\
          & {11} & {100} & {4.46} & {11.00} & {52.83\%} & {51.07\%} & {0.01} & {4.57} & {11.00} & {95.37\%} & {95.20\%} \\
          & {7} & {10} & {4.47} & {7.00} & {52.68\%} & {50.92\%} & {0.01} & {4.70} & {7.00} & {92.59\%} & {92.31\%} \\
    \hline
    \multirow{5}[2]{*}{{(RFE-SVM)}} & {33} & {10} & {4.37} & {33.00} & {94.95\%} & {94.94\%} & {0.01} & {4.50} & {33.00} & {96.75\%} & {96.63\%} \\
          & {22} & {10} & {4.72} & {22.00} & {94.95\%} & {94.94\%} & {0.01} & {4.66} & {22.00} & {96.63\%} & {96.50\%} \\
          & {17} & {100} & {4.48} & {17.00} & {94.95\%} & {94.94\%} & {0.01} & {4.61} & {17.00} & {96.59\%} & {96.46\%} \\
          & {11} & {0.1} & {4.49} & {11.00} & {94.78\%} & {94.76\%} & {0.1} & {4.60} & {11.00} & {96.98\%} & {96.87\%} \\
          & {7} & {0.1} & {4.48} & {7.00} & {94.78\%} & {94.76\%} & {0.1} & {4.83} & {7.00} & {98.09\%} & {98.02\%} \\
    \hline
    (RL-$\ell_1$-M) & -     & 10.00    & 1806.04 & 15.80  & 99.98\% & 99.97\% &    0.1   & 1807.20 & 10.78 & 96.53\% & 96.40\%   \\
    \hline
    (SVM-$\ell_1$) & -     & 10.00    & 3.80  & 17.30  & 99.95\% & 99.95\% & 0.1   & 3.80  & 11.40 & 96.53\% & 96.41\% \\
    \hline
    \end{tabular}%
  }
    \caption{Best average ACC and AUC for the Lepiota dataset.}
  \label{tab:TFCVLepiota}%
  
\end{table}%

The classification performance comparison for the Wdbc dataset is shown in Table~\ref{tab:TFCVWdbc}. In the first case, i.e. 5\% of label noise, 
almost all tested formulations provided nearly the same accuracy and area under the curve with the exception of $\nameM$, which reported the largest classification metrics, selecting only 16.4 features on average. In the second case, the best classifiers are $\nameM$ and $\text{(RL-$\ell_1$-M)}$ reporting 97.54\% ACC and 97.06\% AUC.

The computational results of Mfeat dataset are depicted in Table \ref{tab:TFCVMfeat}. As can be observed, the classification metrics are almost the same for all tested formulations  when 5\% of the dataset contains label noise. 
The one that provided the highest ACC (100\%) and AUC (100\%) is $\text{(FS-SVM)}$. On the other hand, when 5\% of the dataset are SVM outliers, the best classifier is  $\nameM$ reporting 99.90\% ACC and 99.50\% AUC.

Finally, the computational experiments carried out on the  Lepiota dataset are reported in Table \ref{tab:TFCVLepiota}. We can appreciate that the best classification performance  
is obtained by $\nameM$ when 5\% of the dataset contains label noise, with the following metrics: 100\% ACC and 100\% AUC on average. Moreover, this classifier analyzed only 10.7 features on average. In the same manner, when 5\% of the dataset are SVM outliers, the best classifier is $\nameM$ with more than one and a half unit difference in percentage compared to the rest of the models. It reported 99.73\% ACC and 99.72\% AUC. Observe that this classifier selected only 6.4 features on average.

Therefore, in this subsection we have demonstrated  the efficiency of the proposed classifier. We have observed that in the majority of cases, $\nameM$ provided the highest accuracy and area under the curve using less information to classify a new individual than \text{(RL-$\ell_1$-M)}. Note also that the performance of \text{(SVM-$\ell_1$)} is significantly influenced by outliers, providing in several cases the worst results despite selecting many features. Although \text{(FS-SVM)} reported quite good results with label noise, its efficiency is decreased when the dataset contains SVM outliers. (Fisher-SVM) and (RFE-SVM) presents good results for some dataset but their behaviour is not as robust as our model. 

In summary, choosing the parameters which provide the best accuracy in each dataset when contains 5\% of label noise, the average (maximum) percentages of improvement for ACC and AUC of our model with respect the others in the studied datasets are 2.63\% (24.54\%) and 2.94\% (22.78\%), respectively. Similarly,  when the dataset contains 5\% of SVM outliers, the average (maximum) percentages of improvement for ACC and AUC are 2.06\% (14.46\%) and 2.36\% (14.58\%), respectively. Choosing the parameters which provide the best ACC in each dataset, the average (maximum) percentages of improvement for ACC and AUC with respect to each model is depicted in Table \ref{tab:Improvement}. 
The aforementioned results show that $\nameM$ is a very robust classification method improving the existing ones or being among the best in terms of ACC and AUC. Moreover model $\nameM$ involves a reduced number of features in the obtained classifier in contrast to the models (RL-$\ell_1$) and (SVM-$\ell_1$).

\begin{table}[htbp]
  \resizebox{\textwidth}{!}{
  \centering
    \begin{tabular}{l|r|r|r|r}
    \hline
          & \multicolumn{2}{c|}{5\% Label noise} & \multicolumn{2}{c}{5\% SVM outliers} \\
    \hline
          & \multicolumn{1}{p{8.3em}|}{Av. ACC impr. (Max. ACC impr.)} & \multicolumn{1}{p{8.3em}|}{Av. AUC impr.  (Max. AUC impr.)} & \multicolumn{1}{p{8.3em}|}{Av. ACC impr. (Max. ACC impr.)} & \multicolumn{1}{p{8.3em}}{Av. AUC impr. (Max. AUC impr.)} \\
    \hline
     (FS-SVM) & 0.44 (2.81) & 0.32 (3.08) & 2.71(10.13) & 3.70 (14.58) \\
    \hline
    (Fisher-SVM) & 5.50 (23.58) & 5.00 (22.78) & 3.39 (14.46) & 3.15 (11.22) \\
    \hline
    (RFE-SVM)   & 5.49 (8.08) & 7.01 (19.35) & 1.11 (2.54) & 1.26 (2.82) \\
    \hline
    (RL-$\ell_1$-M) & 0.51 (3.88) & 0.75 (4.29) & 0.50 (3.32) & 0.64 (3.64)  \\
    \hline
    (SVM-$\ell_1$) & 1.19 (3.17) & 1.66 (3.21) & 2.58 (5.97) & 3.06 (6.53)  \\
    \hline
    \end{tabular}%
    }
  \caption{Improvement of $\nameM$ with respect to the rest of the models.}
  \label{tab:Improvement}%
  
\end{table}%

\section{Conclusions}
\label{conclusions}
In this paper, we have developed a new model based on support vector machines especially designed for datasets with a large number of features that may contain outliers. We have formulated the model as a mixed-integer problem and proposed an exact strategy for computing the big M parameters of the formulation. Moreover, we have developed a heuristic algorithm to solve it efficiently. We have also validated the heuristic, proving that the obtained upper bound is of high quality.  

Furthermore, we have compared the performance of the proposed classifier with other classifiers based on support vector machines that deal with feature selection or with outlier detection. We have showed the efficiency of our model, whose competitive advantage is that it deals simultaneously with both aspects. Finally, we think that analyzing the proposed model using other $\ell_p$-norms would be interesting for future research.  

\section*{Declaration of interest}

The authors declare that they have no known competing financial interests or personal relationships that could have appeared to influence the work reported in this paper.

\section*{Credit authorship contribution statement}

\textbf{Marta Baldomero-Naranjo:} Conceptualization, Methodology, Software, Writing. \textbf{Luisa I. Martínez Merino:} Conceptualization, Methodology, Software, Writing. \textbf{Antonio M. Rodríguez-Chía:} Conceptualization, Methodology, Software, Writing.

\section*{Acknowledgements}

The authors thank the \textit{Agencia Estatal de
	Investigaci\'on (AEI) and the European Regional Development's funds (ERDF)}: project  MTM2016-74983-C2-2-R, 2014-2020 ERDF Operational Programme and the Department of Economy, Knowledge, Business and University of the Regional Government of Andalusia: projects FEDER-UCA18-106895 and P18-FR-1422, Fundación BBVA: project NetmeetData (Ayudas Fundación BBVA a equipos de investigación científica 2019), and
  \textit{Universidad de C\'{a}diz}: PhD grant UCA/REC01VI/2017 and Programa de Fomento e Impulso de la actividad Investigadora UCA (2018).  The authors would like to thank the anonymous reviewers for their valuable comments and suggestions.

\end{document}